\documentclass[draft]{general}
\usepackage{multirow}
\usepackage{lipsum}
\usepackage{amsfonts}
\usepackage{graphicx}
\usepackage{subcaption}
\usepackage{epstopdf}
\usepackage{algorithm}
\usepackage{algorithmic}
\usepackage{longtable}  
\usepackage{pgfplots}
\usepackage{csvsimple}
\usepackage{tikz}
\usepackage{setstack}
\usepackage{amsthm}
\usepackage{amssymb} 
\usepackage[utf8]{inputenc}
\usepackage{color}
\usepackage{url}
\usepackage{xcolor}
\usepackage{cleveref}
\usepackage{colortbl}
\usepackage{booktabs}
\usepackage{multirow}
\usepackage{cleveref}
\usepgfplotslibrary{fillbetween}
\ifpdf
  \DeclareGraphicsExtensions{.eps,.pdf,.png,.jpg}
\else
  \DeclareGraphicsExtensions{.eps}
\fi

\newcommand{\ex}{\mathbb{E}}
\newcommand{\pr}{\mathbb{P}}
\newcommand{\siga}{\mathcal{F}_{k - \lambda}}
\newcommand{\bound}{M_{k-\lambda}}

\newcommand{\Cinput}[1]{{\IfFileExists{../bib/linear_res.bib}{\input{../#1}}{\input{./#1}}}}
\IfFileExists{../bib/linear_res.bib}{\newcommand{\tpath}{../csvs}}{\newcommand{\tpath}{./csvs}}
\newtheorem{remark}{Remark}

\crefname{hypothesis}{Hypothesis}{Hypotheses}
\newtheorem{lemma}{Lemma}
\crefname{lemma}{Lemma}{Lemmas}
\newtheorem{definition}{Definition}
\crefname{definition}{Definition}{Definitions}
\newtheorem{theorem}{Theorem}
\crefname{theorem}{Theorem}{Theorems}
\newtheorem{corollary}{Corollary}
\crefname{corollary}{Corollary}{Corollaries}
\newtheorem{proposition}{Proposition}
\crefname{proposition}{Proposition}{Propositions}

\newtheorem{assumption}{Assumption}
\crefname{assumption}{Assumption}{Assumptions}
\newtheorem{experiment}{Experiment}
\crefname{experiment}{Experiment}{Experiments}

\definecolor{white}{HTML}{FFFFFF}
\definecolor{red1}{HTML}{FF9292}
\definecolor{red2}{HTML}{FF6161}
\definecolor{red3}{HTML}{FF3434}
\definecolor{red4}{HTML}{FD0000}

\usepackage{amsopn}

\DeclareMathOperator{\row}{row}

\DeclareMathOperator{\col}{col}
\DeclareMathOperator{\rank}{rank}

\begin{document}

\title[Tracking and Stopping]{Solving, Tracking and Stopping Streaming Linear Inverse Problems}
\author[Pritchard and Patel]{Nathaniel Pritchard $^1$ and Vivak Patel $^2$}

\begin{abstract}
In large-scale applications including medical imaging, collocation differential equation solvers, and estimation with differential privacy, the underlying linear inverse problem can be reformulated as a streaming problem.
In theory, the streaming problem can be effectively solved using memory-efficient, exponentially-converging streaming solvers. 
In practice, a streaming solver's effectiveness is undermined if it is stopped before, or well-after, the desired accuracy is achieved. 
In special cases when the underlying linear inverse problem is finite-dimensional, streaming solvers can periodically evaluate the residual norm at a substantial computational cost. 
When the underlying system is infinite dimensional, streaming solver can only access noisy estimates of the residual. 
While such noisy estimates are computationally efficient, they are useful only when their accuracy is known. 
In this work, we rigorously develop a general family of computationally-practical residual estimators and their uncertainty sets for streaming solvers, and we demonstrate the accuracy of our methods on a number of large-scale linear problems. Thus, we further enable the practical use of streaming solvers for important classes of linear inverse problems.
\end{abstract}

\section{Introduction}\label{sec:Intro}
In a myriad of scientific fields, such as medical imaging \cite{birk2014gpu,kakCT}, boundary element analysis \cite{Brebbia1989BoundaryEA}, and differential privacy \cite{updadhyayDP}, large linear inverse problems can be reformulated as streaming problems that can be effectively solved using memory-efficient, exponentially-converging streaming solvers \cite{Gower_2015,clarksonStreaming,patelBlock}. Specifically, the linear problem is reformulated (and approximated) as determining an unknown, $x^* \in \mathbb{R}^n$, using elements of a random sequence $\lbrace (\tilde A_k, \tilde b_k) : k \in \mathbb{N} \rbrace \subset \mathbb{R}^{p \times n} \times \mathbb{R}^p,$\footnote{It is possible to let the copies of $(\tilde A_k, \tilde b_k)$ have different dimensions, but this would require introducing cumbersome notation at this point.} 
satisfying:
\begin{assumption}\label{assum:consis}
   For all $k \in \mathbb{N}$, $\pr (\tilde{A}_k x^* = \tilde b_k) = 1$; and
\end{assumption}
\begin{assumption}\label{assum:independence}
    $\{(\tilde A_k, \tilde b_k): k \in \mathbb{N}\}$ are independent with distribution $D$.
\end{assumption}
\begin{remark}
Some examples of reformulating the linear inverse problem to a streaming model are detailed in \ref{sec:app}. 
\end{remark}
When $n$ and $p$ are large enough that only one observation pair, $(\tilde A_k, \tilde b_k)$, can be stored in RAM, the streaming problem can be solved using a streaming solver---a solver that makes use of $(\tilde A_k, \tilde b_k)$ to update its approximation of $x^*$ and then discards the observation pair. A set of scalable streaming solvers that effectively confronts the complications of this scenario are known as Generalized Randomized Block Kaczamarz (GBRK) methods \cite{patelBlock}. A GBRK method iteratively updates its current approximation by projecting it onto the hyperplanes specified by the observation pair, which allows it to make use of each observation pair and retain a low memory footprint. Importantly, GBRK methods also have (relatively) efficient, geometric convergence rates \cite{Rebrova,Gower_2015,Takac_Ric,patelBlock}.\footnote{If $n$ is of moderate dimension, then more effective solvers can be generated at the cost of more storage \cite{patel2021implicit}.} 

However, a GBRK method's computational efficiency is undermined if it is stopped before, or stopped well-after, the desired solution accuracy is achieved. Thus, a GBRK method's iterates must be carefully tracked and stopped to realize its promise, yet neither tracking nor stopping is straightforward. In the special case where the linear inverse problem is a large linear system, a simple tracking and stopping approach could be to periodically compute the full residual \cite{NEEDELL2014199}, but this has two limitations: (1) an entire period of updates could be wasted before the method is stopped; and (2) in the setting of large linear systems, computing the entire residual is expensive (see numerical experiments in \cref{sec:exp}). In the general case, a na\"{i}ve tracking and stopping method is to use the residual at the current observation pair, but this has its own limitation: a point estimate is only useful if its accuracy is known. In other words, a na\"{i}ve tracking and stopping method based on the residual at the current observation pair is only useful if we can estimate its uncertainty.

In this work, we rigoroulsy develop a general family of computationally-efficient residual estimators and their uncertainty sets. Our family includes the simplest case of using just a single residual from an observation pair, and allows for using a moving window of previously computed, \textit{dependent} residual estimates to reduce the size of the uncertainty set (see \cref{sec:algor}). To analyze this dependent sequence, we develop a novel analytical technique for sub-Exponential distributions (see \cref{sec:conver}), and show that the sub-Exponential model is appropriate for interesting streaming applications (see \ref{sec:app}). From a practical perspective, we demonstrate the scalability of our methodology by tracking and stopping a streaming solver for a streaming collocation problem with one million fixed points (see \cref{sec:exp}). In the context of non-streaming problems, we show that our methodology provides over a 500 times improvement over tracking methods that periodically compute the full residual (see \cref{sec:exp}). Thus, we further enable the practical use of streaming solvers for important classes of linear inverse problems.
\begin{remark}
While related, our current work substantially extends our previous methodology in \cite{pritchard2023towards}, which focused on least squares problems and required the full residual at each update to track a solver's progress.
\end{remark}

\section{Notation} \label{sec:Notation}

Throughout this paper we use, $\mathbb{E}[\cdot]$ to denote an expectation operation, and $\mathbb{P}(\cdot)$ to denote a probability measure.
We let $\mathcal{F}_{k}$ denote the $\sigma$-algebra for $(\tilde A_1, \tilde b_1),\ldots,(\tilde A_k, \tilde b_k )$. Further, we use $\|\cdot\|_2$ to denote a two norm and $\|y\|_B =  \sqrt{\langle y,B y\rangle}$ to be a norm with respect to some symmetric positive definite matrix $B$. We discern between the quantities that we wish to estimate from the estimators by denoting the estimator with a $\hat{\cdot}$ above the quantity it is estimating.

\section{Problem Formulation \& Algorithm}\label{sec:algor}
Recall, we wish to solve (an approximation to) a linear inverse problem by determining an $x^* \in \mathbb{R}^n$ from a stream, $\lbrace (\tilde A_k, \tilde b_k): k \in \mathbb{N} \rbrace \subset \mathbb{R}^{p \times n} \times \mathbb{R}^p$, satisfying 
\cref{assum:consis,assum:independence}; we use a GBRK method to determine $x^*$ from this stream; and we need to efficiently track and stop this GBRK method. In this section, we describe our technique for accomplishing this task. To do this,  we begin by discussing GBRK methods. In \cref{subsection-subexponential}, we present some additional, natural assumptions on the class of distributions for $D$ to make tracking possible. Finally, in \cref{subsection-solver-track-stop}, we introduce our methodology and discuss its salient properties.

\subsection{Generalized Block Randomized Kaczmarz (GBRK)}
Generalized Randomized Block Kaczmarz (GBRK) methods work by iteratively projecting a solution along hyperplanes formed by the row space of independent random observation pairs \cite{Gower_2015,Haddock2022,patel2021implicit}. To be specific, given an initial point $x_0 \in \mathbb{R}^n$ and an inner product on $\mathbb{R}^n$ determined by a symmetric positive definite matrix $B \in \mathbb{R}^{n \times n}$ (usually just the identity), GBRK methods produce iterates recursively by
\begin{equation}\label{eq:g1}
    x_{k+1} = x_k - B^{-1}\tilde A_{k+1}^\top (\tilde A_{k+1}B^{-1}\tilde A_{k+1}^\top)^\dagger (\tilde A_{k+1}x_k - \tilde b_{k+1}).
\end{equation}
GBRKs are known to have geometric convergence rates \cite{Gower_2015,Rebrova,Necora,patelBlock,strohmer}.

\subsection{Sub-Exponential Distribution} \label{subsection-subexponential}

In order to have an uncertainty set for the residual estimate, we need a model for the distribution of the residuals. To understand why, suppose $\tilde A=\tilde b \in \mathbb{R}$ and $\tilde A$ has Pareto distribution with shape parameter of $1$. Then, this pair clearly satisfies \cref{assum:consis}, and $\ex [ |\tilde A x - \tilde b | ] = \infty$ for any $x \neq 1$. In this example, the (absolute) residual at all $x \neq 1$ has arbitrary variability, and, consequently, does not offer any reliable information about the system at any point besides $1$ (even if we took an average over multiple independent copies of $(\tilde A, \tilde b)$). Therefore, to avoid such pathological behavior, we will need to assume some control over the variability of the observations.

We use sub-Exponential distributions as a model for $(\tilde A_k, \tilde b_k)$, which is valid for interesting examples as demonstrated in \ref{sec:app}. A sub-Exponential distribution for a random variable is defined as follows.
\begin{definition}\label{def:SE}
    For a random variable $Y$, with $\ex[Y]= \mu$, $Y-\mu$ follows a sub-Exponential, $\text{SE}(\sigma,\omega)$, distribution with parameters $\sigma$ and $\omega$ if for all $\delta \geq 0$
\begin{equation}
\mathbb{P}\left(|{Y - \mu}| > \delta \right) < 2e^{-\min\left\{\delta^2/(2\sigma^2),\delta/(2\omega) \right\}}.
\end{equation}
    Equivalently, a random variable $Y$ is sub-Exponential, $\text{SE}(\sigma,\omega)$, if 
    \begin{equation}
        \ex [e^{t(Y-\mu)}] \leq e^{\frac{t^2\sigma^2}{2}},
    \end{equation}
    when $|t| < 1/\omega$ \cite{wainwright_2019}.
\end{definition}
These sub-Exponential random variables continue to be sub-Exponential even when scaled by a constant, specifically we have the following lemma, which is a slight modification of the Bernstein inequality in \cite{vershynin_2018}. 
\begin{lemma}\label{lemma:cons_se}
    If given a random variable $Y$  with $\ex[Y] = \mu$  such that $Y - \mu \sim \text{SE}(\sigma,\omega)$ then for constants $c_1,c_2$ where $c_2 \geq c_1 > 0$ it is the case that $c_1(Y -\mu) \sim \text{SE}(c_2\sigma,c_2 \omega)$.
\end{lemma}

The final key property is that all bounded random variables are sub-Exponential, which will be used to show that many example streaming problems are sub-Exponential. Specifically, if we adopt the convention for any positive constant, $c$, $c/0 = \infty$. Then, we have,

\begin{lemma}\cite[Example 2.4]{wainwright_2019}\label{lemma:boundedSE}
    If $Y$ is a random variable with $\ex[Y] = \mu \neq 0$ and $(Y-\mu)/\mu$ takes values in $[y_1,y_2]$, then 
    \begin{equation}
        \frac{Y - \mu}{\mu} \sim \text{SE}\left(\frac{y_2-y_1}{2},0\right).
    \end{equation}
\end{lemma}
With these facts established, we can now state our assumption for the variability of $(\tilde A_k, \tilde b_k)$. 
\begin{assumption}\label{assum:SE}
    There exist $\sigma, \omega \geq 0$ such that, $\forall x\in \mathbb{R}^{n}, k \in \mathbb{N}$,
    \begin{eqnarray}
    &\|\tilde{A}_k x - \tilde b_k\|_2^2 - \ex[\|\tilde{A}_k x - \tilde b_k\|_2^2] \nonumber\\
    &~~~\sim \text{SE}\bigg(\sigma\ex[\|\tilde A_kx - \tilde b_k\|_2^2], {\omega}\ex[\|\tilde A_kx - \tilde b_k\|_2^2]\bigg).
    \end{eqnarray}
\end{assumption}

As detailed in \ref{sec:app}, \Cref{assum:consis,assum:independence,assum:SE} hold in a wide variety of situations including under the instances where the matrices are sketched according to the Johnson-Lindenstrauss transform, in the case when random subsets of a matrix are taken (Randomized Block Kaczmarz), and for a collocation problem evaluated at random points from a hyper-cube.

\subsection{Streaming solver with tracking and stopping} \label{subsection-solver-track-stop}

\Cref{algo:Adaptive-window} describes our method for tracking and stopping streaming solvers applied to streaming reformulations of linear inverse problems whose observations satisfy \cref{assum:consis,assum:independence,assum:SE}. We highlight key aspects of the algorithm below.

\begin{small}
\begin{algorithm}[!htbp]
    \caption{Tracking and Stopping for Least Squares}\label{algo:Adaptive-window}
    \begin{algorithmic}[1]
    \REQUIRE{$B^{-1} \in \mathbb{R}^{n \times n}$ (usually identity matrix), $x_0 \in \mathbb{R}^n$.}
    \REQUIRE{$\lbrace (\tilde A_k, \tilde b_k) \rbrace \subseteq \mathbb{R}^{p \times n} \times \mathbb{R}^p$ satisfying \cref{assum:consis,assum:independence,assum:SE}.}
    \REQUIRE{Moving average window width $\lambda_1 \in \mathbb{N}$.}
    \REQUIRE{$\alpha > 0, \xi_I \in (0,1) , \xi_{II} \in (0,1) , \delta_I \in (0,1)$, $\delta_{II} > 1$, $\eta \geq 1$, $\upsilon > 0$.}
    \STATE $k \leftarrow 0$
    \WHILE{$k == 0$ \OR $\hat \rho_{k} \geq \upsilon$ \OR 
        \begin{eqnarray} \sqrt{\hat \iota_{k}} &\geq \min\Bigg\{\frac{\lambda\eta (1 - \delta_I)^2\upsilon^2}{2\log(1/\xi_I)\sigma^2 \sqrt{\hat \iota_{k}} (1+ \log(\lambda))},\frac{\lambda\eta\upsilon(1-\delta_I)}{2\log(1/\xi_I) \omega }, \nonumber\\ &\quad \quad \quad \quad \frac{\lambda\eta (\delta_{II}-1)^2\upsilon^2}{2\log(1/\xi_{II})\sigma^2 \sqrt{\hat \iota_{k}} (1+ \log(\lambda))},\frac{\lambda\eta\upsilon(\delta_{II} - 1)}{2 \log(1/\xi_{II})\omega } \Bigg\}\nonumber \end{eqnarray}}\label{stopping-condition}
    \STATE \# Iteration $k+1$ \#
    \STATE Receive $\tilde A_{k+1}$ and $\tilde b_{k+1}$
    \STATE $\tilde r_{k+1} \leftarrow \tilde A_{k+1}x_k - \tilde b_{k+1}$
    \IF{$\lambda = 1$}
    		\STATE $\hat \rho_{k+1}, \hat \iota_{k+1} \leftarrow \Vert \tilde r_{k+1} \Vert_2^2, \Vert\tilde r_{k+1} \Vert_2^4$
        \IF{$\|\tilde r_{+1}\|_2^2 > \|\tilde r_{k}\|_2^2$} 
            \STATE $\lambda \leftarrow 2$
        \ENDIF
    \ELSE  
        \STATE $\hat \rho_{k+1}  \leftarrow \sum_{i=k-\lambda +2}^{k+1} \frac{\Vert \tilde r_i \Vert_2^2}{\lambda} $ 
        \STATE $\hat \iota_{k+1} \leftarrow \sum_{i=k-\lambda +2}^{k+1} \frac{\Vert \tilde r_i \Vert_2^4}{\lambda}$
        \IF{$\lambda < \lambda_1$}
            \STATE $\lambda \leftarrow \lambda + 1$
        \ENDIF
    \ENDIF

    \STATE Update the estimated $(1-\alpha)$-interval by computing: $$\hat \rho_{k+1} \pm \Bigg\{\begin{array}{ll}
        \sqrt{2 \log(2/\alpha)  \frac{\sigma^2 \hat \iota_{k+1} (1+\log(\lambda))}{\eta \lambda}} & \text{if }\log(2 / \alpha) \leq \frac{\lambda \sigma^2 (1 + \log(\lambda))}{2\omega^2}\\
        \frac{2\log(2/\alpha)\omega\sqrt{\hat \iota_{k+1}}}{\sqrt{\eta} \lambda} & \text{if } \log(2 / \alpha) > \frac{\lambda \sigma^2 (1 + \log(\lambda))}{2\omega^2}. 
    \end{array} $$
    \label{Credible-Interval}
    
    \STATE $x_{k+1} \leftarrow x_k - B^{-1}\tilde A_{k+1}^{\top} (\tilde A_{k+1} B^{-1} \tilde A_{k+1}^\top)^\dagger \tilde r_{k+1}$

    \STATE $k \leftarrow k + 1$
    \ENDWHILE
    \RETURN{$x_{k}$ \AND estimated $(1-\alpha)$-interval}
    \end{algorithmic}
\end{algorithm}
\end{small}
\subsubsection{Progress Tracking Point Estimate.} 
We track the progress of the solver through the point estimate, 
\begin{equation}
    \hat \rho_k^\lambda  = \sum_{i = k-\lambda +1}^k \frac{\|\tilde A_{i} x_{i-1} - \tilde b_{i}\|_2^2}{\lambda}  = \sum_{i = k-\lambda +1}^k \frac{\|\tilde r_{i}\|_2^2}{\lambda}.
\end{equation} 
Note, $\tilde r_i$ are already computed for the GBRK update, and so our method only requires computing their norm and storing $\lambda$ such scalars. In other words, $\hat \rho_k^\lambda$ requires a marginal computational and memory footprint relative to the overall procedure.
 
When $\lambda = 1$, we see that this point estimate reduces simply to the observed residual at the current iterate. Yet, we have several benefits for considering $\lambda > 1$. 
\begin{enumerate}
\item When $\lambda > 1$, we shrink the uncertainty set as we might expect with independent observations, but with an extra logarithmic term on $\lambda$ because of the statistical dependency between consecutive residual estimates (see \cref{sub-Gaussian:estimator}).
\item The expected value of $\hat \rho_k^\lambda$ (relative to the observation pairs at a given update),
\begin{eqnarray}
    \rho_k^\lambda  &= \sum_{i = k-\lambda +1}^k \frac{\ex [\|\tilde A_{i} x_{i-1} - \tilde b_{i}\|_2^2 | \mathcal{F}_{i-1}]}{\lambda},
\end{eqnarray}
seems to better track the absolute error of the iterate. Indeed, \fref{fig:error:residual} plots the behavior $\rho_k^1$ and $\rho_k^{100}$ against the absolute error of the iterate for a GBRK method applied to a linear system generated from the Phillips matrix in Matrix Depot \cite{matrixdepot}, which shows that the estimator with $\lambda > 1$ corresponds better to the absolute error in comparison to the case when $\lambda = 1$.
\item In a similar context, consider applying a GBRK method to 15 linear systems generated from matrices in Matrix Depot \cite{matrixdepot}. Set a residual threshold to $\Vert Ax_\kappa - b \Vert_2$ where $\kappa = \min\lbrace k : \Vert x_{k} - x^* \Vert_2 <1 \rbrace$. \Fref{fig:stopping:error} plots the number of times $\rho_k^1$ and $\rho_k^{100}$ fall below the residual threshold but $k < \kappa$, which is an imperfect, but informative proxy for determining if the GBRK is stopped too soon. We see that, for nearly half of the systems, $\rho_k^1$ would have stopped early, whereas $\rho_k^\lambda$ does not.
\end{enumerate}

This latter behavior motivates the need for quantifying the uncertainty of $\hat \rho_k^\lambda$ so that we can account for the probabilities of false positives and false negatives.

\begin{figure}
    \begin{subfigure}[t]{.45\textwidth}
        \centering
        \pgfdeclarelayer{background}
\pgfdeclarelayer{foreground}
\pgfsetlayers{background,main,foreground}
\begin{tikzpicture}[baseline]

    \begin{axis}[
        width = \textwidth,
        height = .3\textheight,
        every axis plot/.append style={ultra thick},
        ymode = log,
        xlabel = {$\|x -x^* \|_2^2$},
        label style={font=\small},
        title style={font=\small},
        ylabel style = {yshift = -.25cm},
        xlabel style = {yshift = .2cm},  
        tick label style={font=\small},
        ytick pos = left,
        xtick pos = bottom]
        \begin{pgfonlayer}{background}
            \addplot[red!60!black, only marks, on layer = background, draw opacity = 0, fill opacity = .3] table [x = Error, y = True_Res, col sep=comma]{\tpath/phillips.csv};
            \addplot[blue, only marks, on layer = background, draw opacity = 0, fill opacity = .3] table [x = Error, y = Mov_True_Res, col sep=comma]{\tpath/phillips.csv};
        \end{pgfonlayer}
         
      \legend{$\lambda =1$, $\lambda >1$}
    \end{axis}
    \draw[white] (0,-1.9) rectangle (4.5,-.1);
\end{tikzpicture}
        \caption{Plot of the $\rho_k^\lambda$ when $\lambda = 1$ (red) and $\lambda > 1$ (blue) at different values of error for a Phillips matrix.}\label{fig:error:residual}
    \end{subfigure}
    \hspace{1em}
    \begin{subfigure}[t]{.45\textwidth}
        \centering
         \begin{tikzpicture}[baseline]
\begin{axis}[
    width = \textwidth,
    height = .3\textheight,
    ybar,
    ylabel=Early Stops,
    ylabel near ticks,
    xtick={1,2,3,4,5,6,7,8,9,10,11,12,13,14,15},
    xticklabels={frank,moler,hankel, chow,circul,parter, kms, dingdong,grcar,hadamard, pei, rohess, sampling, triw, wilkinson},
    xticklabel style={rotate=90},
    bar width=0.1cm,
    xtick pos = bottom,
    ytick pos = left
 ]
    \addplot[bar shift = 0cm, fill=red!60!black, color = red!60!black] coordinates{(1, 21) (2,18)(3, 10) (4, 6) (5, 4) (6, 6)(7,3) (8, 0)   (9, 0) (10, 0)  (11, 0) (12,0) (13,0) (14, 0) (15, 0)};
    \addplot[fill=blue, color = blue, bar shift = .2cm] coordinates{ (1, 0)  (2,0) (3, 0) (4, 0) (5, 0) (6, 0) (7,0) (8, 0) (9, 0) (10, 0)  (11, 0) (12,0) (13,0) (14, 0) (15, 0)}; 
 \legend{$\lambda =1$, $\lambda >1$}
\end{axis}
\end{tikzpicture}
        \caption{Graph showing the number of iterates with a $\rho_k^\lambda$ less than the $\rho_{\kappa}^\lambda$  where $\kappa=\min\lbrace k : \Vert x_{k} - x^* \Vert_2 <1\rbrace$ for $\lambda = 1$ (red) and $\lambda > 1$ (blue).} \label{fig:stopping:error}
    \end{subfigure}
\end{figure}

\subsubsection{Progress Tracking Uncertainty Set Estimate.}  \label{subsub-uq}
In line 18, a $(1-\alpha)$ uncertainty set for $\hat \rho_{k}^\lambda$ is computed, where $\alpha \in (0,1)$ is given by the user (and reflect the user's risk tolerance). This uncertainty set is an estimate of the true uncertainty set for $\hat \rho_k^\lambda$ derived in \cref{coll:credible-interval}. The true uncertainty set satisfies either 
\begin{equation} \label{eq:credible-interval-def-1}
\pr\left( \rho_k^\lambda \in \hat \rho_k^\lambda \pm \sqrt{2 \log(2/\alpha)  \frac{\sigma^2 \hat \iota_{k} (1+\log(\lambda))}{\lambda}} \right) \geq 1-\alpha, 
\end{equation}
when $\log(2 / \alpha) \leq \frac{\lambda \sigma^2 (1 + \log(\lambda))}{2\omega^2}$; or
\begin{equation} \label{eq:credible-interval-def-2}
\pr\left( \rho_k^\lambda \in \hat \rho_k^\lambda \pm \frac{2\log(2/\alpha)\omega\sqrt{\hat \iota_{k}}}{ \lambda} \right) \geq 1-\alpha,
\end{equation}
otherwise.

Compared to \eref{eq:credible-interval-def-1} and \eref{eq:credible-interval-def-2}, line 18 has an additional $\eta$ term. This term is included to adjust for conservativeness of the theoretical intervals. Choices of $\eta$ for specific distributions can be found in \ref{ex:diff_priv}.

\subsubsection{Stopping Criterion} \label{subsub-stop}
After, the completion of the uncertainty quantification stage, the solution is updated in lines 19-20 in accordance with \eref{eq:g1}. Then on line 2 the stopping criterion is checked.  Note, in a deterministic setting stopping is simple: stop when the tracking value falls below some threshold $\upsilon$. When randomness is introduced, this criterion can incur two different errors.
The first can be viewed as stopping too late, and it occurs when the tracking parameter value, $\rho_k^\lambda \leq \delta_I \upsilon$, while $\hat \rho_k^\lambda > \upsilon$, where $\delta_I$ is a user defined parameter that permits the specification of where the gap between $\hat \rho_k^\lambda$ and $\rho_k^\lambda$ is large enough to be considered problematic. By using the condition on $\hat{\iota}^\lambda_k$ (defined on line 13), we approximately control the probability of this error at $\xi_{I}$ (see \cref{corr:stopping} and \Cref{subsec:est_stop}).

The second error type can be viewed as stopping too early, and it occurs when the tracking parameter value, $\rho_k^\lambda \geq \delta_{II} \upsilon$, while $\hat \rho_k^\lambda < \upsilon$ (c.f., \fref{fig:stopping:error}), Where $\delta_{II}$ is a user defined parameter that permits the specification of where the gap between $\hat \rho_k^\lambda$ and $\rho_k^\lambda$ is great enough to be considered problematic. By choosing the right stopping criterion we can then control the probability of this error at  $\xi_{II}$.

As we will show in \cref{corr:stopping}, the criterion dependent on $\hat \iota_k^\lambda$ in line 2 of the algorithm accomplishes this task.  Once this condition is satisfied, then the deviations between $\rho_k^\lambda$ and $\hat \rho_k^\lambda$ are reasonably well controlled; thus, it is safe to stop when $\hat \rho_k^\lambda < \upsilon$ (see \cref{corr:stopping} and \Cref{subsec:est_stop}). 

\subsubsection{Moving Average Logistics}
During the early phase of the algorithm, a GBRK usually experiences a rapid convergence to a region of the solution. During this phase, the GBRK's estimated residual reflects this rapid convergence and warrants a small moving average window. At later iterations, the convergence slows down and most of the variability in the residuals comes from randomness. During this phase, the moving average windows should be larger. Lines 6-10 and 14-15 reflect this behavior by starting the moving window at $\lambda = 1$, and increasing it to $\lambda$ once the phase change is detected. The difference between these phases is determined to be the iteration $k'$ where $\|\tilde r_{k'}\|_2^2 > \|\tilde r_{k'-1}\|_2^2$, as seen on line 8, which yields good empirical behavior.

\section{Consistency of Estimators and Uncertainty Sets}\label{sec:conver}
Core to the establishment of the theoretical validity of \cref{algo:Adaptive-window} is proving the consistency of the estimator of, and reliability of the uncertainty sets for, $\hat \rho_k^\lambda$. 
To accomplish these two tasks, it is first necessary to show that the general form of \eref{eq:g1} converges in all moments, which is novel in comparison to previous analyses of GBRKs. 
From there, we can then combine this convergence result with Chernoff bounds to derive the distribution around $\hat \rho_k^\lambda$, which can be used to show $\hat \rho_k^\lambda$'s consistency for $\rho_k^\lambda$. With the consistency of $\hat \rho_k^\lambda$ established, we use the distribution of $\hat \rho_k^\lambda$ to derive its uncertainty set. Unfortunately, the uncertainty set will rely on an uncomputable quantity, so we close the section by showing that the uncomputable quantity can be estimated by $\hat \iota_k^\lambda$ with a reasonable relative error.
 
\subsection{Convergence of the Residuals' Moments}\label{subsec:allmom}

Our first goal is to prove that all the moments of the residual will converge to zero. To achieve this goal, we will show that the moments of the absolute error converge to zero. To this end, we will transform the iteration update, \eref{eq:g1}, into a more amenable form. As we will see, this more amenable form shows that the updates are a sequence of orthogonal projections. Using these orthogonal projections and \cref{Thm:Patel:Im:4.1}, we will show that, at a random iteration, a sufficient, random geometric reduction in the error will occur. Our final step will be to control the random iteration and the random reduction in the error. Once these pieces are in place, we will be able to conclude that the moments of the absolute error and, hence, the residual, decay (geometrically) to zero.

\paragraph{Transformation of Variables}
To avoid unnecessary considerations about inner products, we will begin with a transformation of the variables by a symmetric square root of $B$. In other words, \eref{eq:g1} becomes
\begin{small}
\begin{equation} \label{eq:transformed-update}
B^{1/2} x_{k+1} = B^{1/2} x_k - B^{-1/2}\tilde A_{k+1}^\top(\tilde A_{k+1}B^{-1}\tilde A_{k+1}^\top )^\dagger (\tilde A_{k+1} x_k - \tilde b_{k+1}).
\end{equation}
\end{small}
To simplify this relationship further, it will be useful to introduce several important spaces. Let
\begin{small}
\begin{eqnarray}
\mathcal{H} &=\lbrace x \in \mathbb{R}^n, \forall k \in \mathbb{N}: \mathbb{P}(\tilde A_k x = \tilde b_k) = 1 \rbrace, \\
\mathcal{N} &=\lbrace x \in \mathbb{R}^n, \forall k \in \mathbb{N}: \mathbb{P}(\tilde A_kB^{-1/2}x = 0) = 1 \rbrace,~\text{and} \\
\mathcal{R} &= \mathcal{N}^\perp.
\end{eqnarray}
\end{small}
Under \cref{assum:consis}, $\mathcal{H} \neq \emptyset$ and denotes the set of all solutions to the linear inverse problem. Moreover, $\mathcal{N}$ represents the null space of the linear problem, which can be equivalently written as $\mathcal{N} = \mathrm{null}( B^{-1/2} \ex [ \tilde A_k^\top \tilde A_k ] B^{-1/2})$. From this characterization, $\mathcal{R}$ represents the row space of the linear problem and, equivalently, $\mathcal{R} = \row( B^{-1/2} \ex [\tilde A_k^\top \tilde A_k ] B^{-1/2})$.

Using these spaces, let $x^*$ be the orthogonal projection of $x_0$ onto the set $\mathcal{H}$, and let $\beta_k = B^{1/2} (x_{k} - x^*)$. Then, \eref{eq:transformed-update} simplifies to
\begin{small}
\begin{equation} \label{eq:transformed-error-iteration}
\beta_{k+1} = \beta_k - B^{-1/2} \tilde A_{k+1}^\top (\tilde A_{k+1}B^{-1}\tilde A_{k+1}^\top )^\dagger \tilde A_{k+1} B^{-1/2}\beta_k.
\end{equation}
\end{small}

From \eref{eq:transformed-error-iteration}, we observe that $B^{-1/2} \tilde A_{k+1}^\top  (\tilde A_{k+1}B^{-1}\tilde A_{k+1}^\top)^\dagger \tilde A_{k+1}B^{-1/2}$ is an orthogonal projection onto $\row(\tilde A_{k+1} B^{-1/2})$. As a result, we have an observation and a useful simplification. First,
\begin{small}
\begin{lemma} \label{lemma:transformed-errors-are-in-rowspace}
$\lbrace \beta_k : k+1\in\mathbb{N} \rbrace \subset \mathcal{R}$.
\end{lemma}
\end{small}
\begin{proof}
 By construction, $\mathcal{N} \subset \mathrm{null}(\tilde A_k B^{-1/2})$ for all $k$ with probability one. Hence, $\row(\tilde A_k B^{-1/2}) \perp \mathcal{N}$ with probability one. Letting $\mathcal{P}_\mathcal{N}$ denote the orthogonal projection operator onto $\mathcal{N}$, if $\mathcal{P}_\mathcal{N} \beta_{k} = 0$ then $\mathcal{P}_\mathcal{N} \beta_{k+1} = 0$ by \eref{eq:transformed-error-iteration}. Noting that from our definition of $x^*$, $\beta_{0} = B^{-1/2}(x_0-x^*)$ can only be in $\mathcal{N}$ when $x_0 = x^*$. By construction, since $\beta_{0} \not\in \mathcal{N} \Rightarrow \beta_{0} \in \mathcal{N}^{\perp}$; thus, $\mathcal{P}_{\mathcal{N}} \beta_{0} = 0$. The result follows by induction.
\end{proof}
Second, if we let $Q_{k+1}$ be a matrix with orthonormal columns that form a basis for $\row(\tilde A_{k+1} B^{-1/2})$, then \eref{eq:transformed-error-iteration} becomes
\begin{small}
\begin{equation} \label{eq:transformed-error-iter-ortho}
\beta_{k+1} = \beta_k - Q_{k+1} Q_{k+1}^\top \beta_k.
\end{equation}
\end{small}
\paragraph{Geometric Reduction in Error}

Now, let $\tau_0 = 0$ and let $\tau_1$ be the first iteration such that 
\begin{equation}
 \col(Q_{1})+ \col(Q_{2}) + \cdots + \col(Q_{\tau_1})  = \mathcal{R},
\end{equation}
otherwise let $\tau_1$ be infinite. When $\tau_1$ is finite, we can use the following extension of Meany's lemma proposed in \cite[Theorem 4.1]{patel2021implicit} about the convergence of sequences of orthogonal projections.
\begin{theorem}[\cite{patel2021implicit}, Theorem 4.1] \label{Thm:Patel:Im:4.1}
    Let $q_1,\dots,q_k$ be unit vectors in $\mathbb{R}^{n}$ for some $k \in \mathbb{N}$. Let $\mathcal{V} = \text{span}[q_1,\dots,q_k]$. Let $\mathcal{W}$ denote all matrices $W$, where the columns of $W$ are the vectors $\{w_1,\dots, w_k\} \subset  \{q_1,\dots,q_k\}$ that are a maximal linearly independent subset. Then 
    \begin{equation}
        \sup_{y \in \mathcal{V}, \|y\|_2 = 1} \|Qy\|_2 \leq \sqrt{1 - \min_{W\in \mathcal{W}}\text{det}(W^\top W)},
    \end{equation}
    where $Q = (I - q_kq_k^\top) \cdots(I - q_1q_1^\top)$.
\end{theorem}

By combining \cref{Thm:Patel:Im:4.1} with our observation that \eref{eq:transformed-error-iteration} is a sequence of orthogonal projections, we obtain the following lemma.
\begin{lemma}\label{lemma:gamma-between-0-1}
Let $x_0 \in \mathbb{R}^n$. Let $\lbrace x_k \rbrace$ be generated according to \eref{eq:g1} for $(\tilde A_k, \tilde b_k)$ from a distribution $D$ satisfying \cref{assum:consis}. Let $x^*$ be the orthogonal projection of $x_0$ onto $\mathcal{H}$. On the event, $\lbrace \tau_1 < \infty \rbrace$, there exists a $\gamma_1\in (0,1)$ that is a function of $\lbrace Q_1,\ldots,Q_{\tau_1}\rbrace$ such that
\begin{small}
\begin{equation}
\Vert x_{\tau_1} - x^* \Vert_B \leq \gamma_1 \Vert x_0 - x^* \Vert_B.
\end{equation}
\end{small}
\end{lemma}
\begin{proof}
By our definition of $\beta_k$, we need only prove that $\exists \gamma_1 \in (0,1)$ such that $\Vert \beta_{\tau_1} \Vert_2 \leq \gamma_1 \Vert \beta_0 \Vert_2$. To prove this, let $q_{k,1},\ldots,q_{k,p}$ denote the columns of $Q_k$. Then, by \eref{eq:transformed-error-iter-ortho},
\begin{small}
\begin{equation}
\beta_{\tau_1} = \left[ \prod_{k=1}^{\tau_1} \left( \prod_{j=1}^p (I - q_{k,j} q_{k,j}^\top ) \right) \right] \beta_0
\end{equation}
\end{small}
Since $\beta_0 \in \mathcal{R}$ by \cref{lemma:transformed-errors-are-in-rowspace}. \cref{Thm:Patel:Im:4.1} implies that there exists a $\gamma_1 \in (0,1)$ that is a function of $\lbrace q_{1,1}, q_{1,2},\ldots,q_{\tau_1,p-1}, q_{\tau_1,p} \rbrace$ such that $\Vert \beta_{\tau_1} \Vert_2 \leq \gamma_1 \Vert \beta_0 \Vert_2$.
\end{proof}

We need not stop at $\tau_1$. In fact, we can iterate on this argument in the following manner. When $\lbrace \tau_{\ell} < \infty \rbrace$, define $\tau_{\ell+1}$ to be the first iteration after $\tau_{\ell}$ such that
\begin{small}
\begin{equation}
 \col(Q_{\tau_{\ell} + 1}) + \col( Q_{\tau_{\ell}+2} ) + \cdots + \col(Q_{\tau_{\ell+1}})  = \mathcal{R},
\end{equation}
\end{small}
otherwise let $\tau_{\ell+1}$ be infinite. Then, we have the following straightforward corollary.
\begin{corollary} \label{corollary:geometric-convergence}
Let $x_0 \in \mathbb{R}^n$. Let $\lbrace x_k \rbrace$ be generated according to (\ref{eq:g1}) for $(\tilde A_k, \tilde b_k)$ from a distribution $D$ satisfying \cref{assum:consis}. Let $x^*$ be the orthogonal projection of $x_0$ onto $\mathcal{H}$. On the event, $ \cap_{\ell = 1}^L \lbrace \tau_\ell < \infty \rbrace$, there exists $\gamma_\ell\in (0,1)$ that is a function of $\lbrace Q_{\tau_{\ell-1}+1},\ldots,Q_{\tau_{\ell}}\rbrace$ for $\ell = 1,\ldots,L$, such that
\begin{small}
\begin{equation}
\Vert x_{\tau_L} - x^* \Vert_B \leq \left( \prod_{\ell=1}^L \gamma_\ell \right) \Vert x_0 - x^* \Vert_B.
\end{equation}
\end{small}
\end{corollary}

\paragraph{Control of the Random Rate and Random Iteration}
Of course, \cref{corollary:geometric-convergence} does not imply that the absolute error converges to zero. In fact, \cref{corollary:geometric-convergence} has two points of failure. First, it may happen that $\gamma_{\ell} \to 1$ as $\ell \to \infty$ with some nonzero probability; that is, we have no control over the random rate of convergence. This issue is addressed by the following result, which relies on the independence of $\lbrace (\tilde A_k, \tilde b_k) \rbrace$.

To obtain this results we rely on the following key theorem from \cite[Theorem 4.1.3]{durrett2013probability}. 
\begin{theorem}[\cite{durrett2013probability}, Theorem 4.1.3]\label{Thm:durrett:4.1.3}
    Let $Z_1,Z_2,\dots$ be $i.i.d.$ random variables, $\mathcal{Z}_n = \sigma(Z_1,\dots,Z_k)$ and $\tau$ be a stopping time with $\mathbb{P}(\tau < \infty) > 0$. Conditioned on $\left\{\tau < \infty\right\}$, $\{Z_{\tau+k}, k > 1\}$ is independent of $\mathcal{Z}_n$ and has the same distribution as $Z_1,\dots,Z_k$.
\end{theorem}
We now can use \cref{Thm:durrett:4.1.3} to show that $\gamma_\ell \to 1$ does not occur.
\begin{lemma} \label{lemma:convergence-rate-iid}
Let $x_0 \in \mathbb{R}^n$. Let $\lbrace x_k \rbrace$ be generated according to (\ref{eq:g1}) for $\lbrace (\tilde A_k, \tilde b_k) \rbrace$ from a distribution $D$ satisfying \cref{assum:consis,assum:independence}. Let $x^*$ be the orthogonal projection of $x_0$ onto $\mathcal{H}$.
Then, whenever they exist, $\lbrace \tau_{\ell} - \tau_{\ell-1}: \ell \in \mathbb{N} \rbrace$ are independent and identically distributed; and $\lbrace \gamma_{\ell}: \ell \in \mathbb{N} \rbrace$ are independent and identically distributed.
\end{lemma}
\begin{proof}
When $\tau_{\ell}$ is finite, by \cref{Thm:durrett:4.1.3}, $\lbrace Q_{\tau_{\ell} +1}, \ldots, Q_{\tau_{\ell}+k} \rbrace$ given $\tau_{\ell}$ are independent of $\lbrace Q_{1},\ldots,Q_{\tau_{\ell}} \rbrace$ and are identically distributed to $\lbrace Q_1,\ldots,Q_k \rbrace$  for all $k$. Therefore, $\tau_{\ell} - \tau_{\ell-1}$ are independent and identically distributed, as are $\gamma_{\ell}$. 
\end{proof}

The second point of failure in \cref{corollary:geometric-convergence}, as alluded to in \cref{lemma:convergence-rate-iid}, is the existence of $\lbrace \tau_{\ell} \rbrace$. Specifically, on $\lbrace \tau_{\ell} = \infty \rbrace$, \cref{corollary:geometric-convergence} will no longer supply a rate of improvement in the absolute error. Therefore, we must show that $\lbrace \tau_{\ell} < \infty \rbrace$ occurs with probability one, which is the content of the next result.

\begin{lemma} \label{lemma:iteration-control}
Let $x_0 \in \mathbb{R}^n$. Let $\lbrace x_k \rbrace$ be generated according to \eref{eq:g1} for $\lbrace (\tilde A_k, \tilde b_k) \rbrace$ from streaming distribution $D$, satisfying \cref{assum:consis,assum:independence,assum:SE}. Let $x^*$ be the orthogonal projection of $x_0$ onto $\mathcal{H}$.

Then $\mathbb{P}( \tau_\ell < \infty) = 1$ for every $\ell \in \mathbb{N}$. Moreover, $\exists \pi \in (0,1]$ such that, for all $\ell \in \mathbb{N}$ and $k \geq \rank(\mathcal{R})$,
\begin{small}
\begin{equation}\label{lemma:itcontrol:eq:negbinom}
\mathbb{P}( \tau_{\ell} - \tau_{\ell-1} = k ) \leq {k-1\choose \rank(\mathcal{R}) - 1}(1-\pi)^{k - \rank(\mathcal{R})}\pi^{\rank(\mathcal{R})}.
\end{equation}
\end{small}
\end{lemma}
\begin{proof}
Given that $\lbrace Q_k : k \in \mathbb{N} \rbrace$ are independent and identically distributed, we will show that the probability that $\col(Q_1) + \cdots + \col(Q_{k+1})$ grows in dimension relative to $\col(Q_1) + \cdots + \col(Q_k)$, when $\dim(\col(Q_1) + \cdots + \col(Q_k) ) < \rank(\mathcal{R})$ is some $\pi \in (0,1]$. As a result, the probability that the dimension increases $\rank(\mathcal{R})$ times in $k$ iterations (with $k \geq \rank(\mathcal{R})$) is dominated by a negative binomial distribution. In other words, for $k \geq \rank(\mathcal{R})$,
\begin{small}
\begin{equation}
\mathbb{P}\left( \tau_1 = k \right) \leq {k - 1 \choose \rank(\mathcal{R}) - 1}(1 - \pi)^{k - \rank(\mathcal{R})}\pi^{\rank(\mathcal{R})}.
\end{equation}
\end{small}
This implies $\tau_1$ is finite with probability one, and the result for $\tau_{\ell} - \tau_{\ell-1}$ follows by \cref{lemma:convergence-rate-iid}.

Thus, it only remains to show that the probability that the dimension grows is bounded from below by $\pi' \in (0,1]$. To do so, we need only show that $\exists \pi' > 0$ such that for any $z \in \mathcal{R}$, $\mathbb{P}( \Vert Q_1^\top z \Vert_2^2 > 0 ) \geq \pi'$. Since $\pi'$ represents a lower bound on the probability that any vector in $\mathcal{R}$ falls in the column space of a $Q_i$, it implies if we have some vector $z' \in \mathcal{R}$ but $z' \notin \col(Q_1) + \cdots + \col(Q_k)$ then we know $z' \in \col(Q_{k+1})$ with probability at least $\pi'$, which implies that $\dim(\col(Q_1) + \cdots + \col(Q_{k+1}) )$ will increase with at least the same probability. If we now recall that $\col(Q_1) = \row(\tilde A_1 B^{-1/2} ) \subset \mathcal{R}$, then we can note that $\mathbb{P}( \Vert Q_1^\top z \Vert_2^2 > 0 )=\mathbb{P}( \Vert \tilde A_k B^{-1/2} z \Vert_2^2 > 0 )$, and thus if we observe enough stream blocks, the sum of the spaces will span all of $\mathcal{R}$.
We prove the existence of $\pi'$ by noting from \cref{def:SE} that for any $\delta \in (0,1)$,
\begin{small}
\begin{eqnarray}
\mathbb{P}\bigg( \Vert\tilde A_k B^{-1/2} z \Vert_2^2 &> 0 \bigg) 
\geq  \mathbb{P}\left( \Vert\tilde A_k B^{-1/2} z \Vert_2^2 \geq (1 - \delta) \ex [ \Vert\tilde A_k B^{-1/2} z \Vert_2^2 ] \right) \label{eq:sub-exp}\\
&\geq  1 - \mathbb{P}\left( \frac{\Vert\tilde A_k B^{-1/2} z \Vert_2^2 - \ex [ \Vert\tilde A_k B^{-1/2} z \Vert_2^2 ]}{\ex [ \Vert\tilde A_k B^{-1/2} z \Vert_2^2 ]} < - \delta\right) \\
&\geq 1 - \exp( - \min\lbrace \delta^2/(2\sigma^2), \delta/(2\omega) \rbrace),
\end{eqnarray}
\end{small} 
where (26) comes from \cref{assum:SE}.

Note, for any $\delta \in (0,1)$, the last term is strictly larger than $0$. Hence, we can fix a $\delta \in (0,1)$ and set the resulting value to $\pi'$.
\end{proof}

\paragraph{Convergence of the Moments}
We now put these pieces together to conclude as follows.
\begin{theorem}\label{thm:convergence-of-moments}
Let $x_0 \in \mathbb{R}^n$. Let $\lbrace x_k \rbrace$ be generated according to \eref{eq:g1} for $\lbrace (\tilde A_k, \tilde b_k) \rbrace$ from distribution $D$ and satisfying \cref{assum:consis,assum:independence,assum:SE}. Let $x^*$ be the orthogonal projection of $x_0$ onto $\mathcal{H}$.
Then, for any $d \in \mathbb{N}$, $\mathbb{E}[ \Vert \tilde A_{k+1} x_k - \tilde b_{k+1} \Vert_2^d]  \to 0$ and $\mathbb{E}[ \Vert x_k - x^* \Vert_B^d ] \to 0$ as $k \to \infty$. Additionally,
for any $\ell \in \mathbb{N}$, we have 
\begin{small}
\begin{equation}
    \mathbb{E}\left[ \Vert x_{\tau_{\ell}} - x^* \Vert_B^d \right] \leq \mathbb{E}\left[\gamma_1^d \right]^\ell \Vert x_0 - x^* \Vert_B^d.
\end{equation}
\end{small}
\end{theorem}
\begin{proof}
It is enough to show that $\mathbb{E}[ \Vert x_k - x^* \Vert_B^d ] \to 0$ as $k \to \infty$. By \eref{eq:transformed-error-iter-ortho}, the absolute error is a non-increasing sequence. Thus, we need only show that a subsequence converges to zero. By \cref{corollary:geometric-convergence,lemma:convergence-rate-iid,lemma:iteration-control},
\begin{small}
    \begin{equation}
        \mathbb{E}\left[ \Vert x_{\tau_{\ell}} - x^* \Vert_B^d \right] \leq \mathbb{E}\left[\gamma_1^d \right]^\ell \Vert x_0 - x^* \Vert_B^d,
    \end{equation}
\end{small}

for all $\ell \in \mathbb{N}$, where $\mathbb{E}[\gamma_1^d] < 1$. Therefore, as $\ell \to \infty$, the conclusion follows.
\end{proof}

\subsection{Distribution and Consistency of Estimators} \label{subsec:cons}
Using the convergence of the moments established in \cref{thm:convergence-of-moments}, we now determine the distribution of $\hat \rho_k^\lambda$. Determining the distribution of $\hat \rho_k^\lambda$ would be easy if the terms in the sum of the norm squared of the residuals,$\|\tilde r_{i}\|_2^2 = \|\tilde A_{i} x_{i-1} - \tilde b_{i}\|_2^2$, were independent, as \cref{assum:SE} would imply that $ \hat \rho_k^\lambda$ is sub-Exponential. Unfortunately, the terms composing $\hat \rho_k^\lambda$ are dependent. To handle this dependence, we innovate the following proof and conclude, $\hat \rho_k^\lambda$ is sub-Exponentially distributed with a variance that is only a logarithmic term worse than the independent case.

Since we are concerned with the consistency of the estimators we need to first understand the distribution conditioning on $\siga$, the $\sigma$-algebra generated by $\{(\tilde A_j, \tilde b_j): j= 1, \cdots k-\lambda\}$, of $\|\tilde r_i\|_2^2 - \ex[\|\tilde r_i\|_2^2 | \mathcal{F}_{i-1} ] | \mathcal{F}_{k-\lambda}$ for any $i > k - \lambda$. 
\begin{lemma}\label{lemma:abs_diff}
  Suppose the setting of \cref{thm:convergence-of-moments} holds.
Then, for any $i > k - \lambda$, 

    \begin{equation}
        \|\tilde r_i\|_2^2 - \ex[\|\tilde r_i\|_2^2 | \mathcal{F}_{i-1} ] \Big\vert \siga  \sim \textbf{SE}(\sigma \bound, \omega \bound),
    \end{equation}

    where $\bound = \ex [\Vert \tilde A_{k - \lambda} B^{-1/2} \Vert_2^2] \Vert x_{k - \lambda-1} - x^* \Vert_B^2$. 
\end{lemma}
\begin{proof}
By \eref{eq:transformed-error-iter-ortho},
\begin{equation}
    \ex [ \Vert \tilde r_{i} \Vert_2^2 | \mathcal{F}_{i-1}] \leq \ex [\Vert \tilde A _{i}B^{-1/2} \Vert_2^2]\Vert x_{i-1} - x^* \Vert_B^2 \leq \bound \label{rtoM}
\end{equation}

Using \eref{rtoM} and \cref{assum:SE}, for any $\delta > 0$,
\begin{eqnarray}
&\mathbb{P}\left( \left\vert \|\tilde r_i\|_2^2 - \ex[\|\tilde r_i \|_2^2 | \mathcal{F}_{i-1} ] \right\vert > \delta \Big\vert \mathcal{F}_{k-\lambda} \right) \nonumber \\
&< 2\exp\left(- \min \left\lbrace \frac{\delta^2}{2\sigma^2 \ex[\|\tilde r_i \|_2^2 | \mathcal{F}_{i-1}]^2 }, \frac{\delta}{2\omega \ex[\|\tilde r_i \|_2^2 | \mathcal{F}_{i-1}]} \right\rbrace\right) \\
&< 2\exp\left(- \min \left\lbrace \frac{\delta^2}{2\sigma^2 \bound^2 }, \frac{\delta}{2\omega \bound } \right\rbrace\right)
\end{eqnarray}

Taking expectations of the probabilities with respect to $\siga$ in the previous display equation and noting that $\bound$ is measurable with respect to $\siga$, we conclude.  
\end{proof}

Using this lemma, we now characterize the distribution of $\hat \rho_k^\lambda - \rho_k^\lambda$ to be sub-Exponential. 
\begin{theorem}\label{sub-Gaussian:estimator}
    Suppose the setting of \cref{thm:convergence-of-moments} holds. Then, 
    \begin{small}
        \begin{equation}
            \hat \rho_{k}^\lambda - \rho_{k}^\lambda \Big| \mathcal{F}_{k-\lambda} \sim \mathbf{SE}\left(\sigma \bound\sqrt{\frac{(1 + \log(\lambda))}{ \lambda }},\frac{\omega \bound}{\lambda} \right).
        \end{equation} 
    \end{small}

\end{theorem}
\begin{proof}
    Using \cref{lemma:abs_diff} and induction, we will prove, for $|t| \leq \lambda/(\omega \bound)$, 
    \begin{small}
    \begin{eqnarray}\label{eq:chernoff}
    \ex& \left[ e^{ t (\hat \rho_k^\lambda - \rho_k^\lambda )} \bigg| \siga \right]\\ 
    &=\mathbb{E}\left[ \left. \prod_{i=k-\lambda+1}^k \exp\left\lbrace\frac{t}{\lambda} \left( \Vert \tilde r_i \Vert_2^2 - \ex\left[\Vert \tilde r_i \Vert_2^2  | \mathcal{F}_{i-1} \right]\right) \right\rbrace  \right\vert \siga \right] \\
    & \leq \exp\left( \frac{t^2 \sigma^2 \bound^2}{2 \lambda } \sum_{j=1}^\lambda \frac{1}{j} \right).
    \end{eqnarray}
    \end{small}
    We can then use a logarithm to bound the summation, which yields the conclusion.
    
    The base case of $\lambda = 1$ follows trivially from \cref{lemma:abs_diff}. Now assume that the result holds up to $\lambda - 1$. Then,
    \begin{small}
    \begin{eqnarray}
    &\mathbb{E}\left[ \left. \prod_{i=k - \lambda+1}^k \exp\left\lbrace\frac{t}{\lambda} \left( \Vert \tilde r_i \Vert_2^2 - \ex\left[\Vert \tilde r_i \Vert_2^2 | \mathcal{F}_{i-1} \right] \right) \right\rbrace  \right\vert \siga \right] \\
    &= \mathbb{E}\left[ \mathbb{E}\left[\left. \left. \prod_{i=k - \lambda+1}^{k} \exp\left\lbrace\frac{t}{\lambda} \left( \Vert \tilde r_i \Vert_2^2 - \ex\left[\Vert \tilde r_i \Vert_2^2 | \mathcal{F}_{i-1} \right] \right) \right\rbrace \right\vert \mathcal{F}_{k-1} \right] \right\vert \siga \right] \\
    &=\mathbb{E}\left[  \mathbb{E} \left[ \left. \exp\left\lbrace\frac{t}{\lambda} \left( \Vert \tilde r_k \Vert_2^2 - \ex\left[\Vert \tilde r_k \Vert_2^2  | \mathcal{F}_{k-1} \right]\right) \right\rbrace \right\vert \mathcal{F}_{k-1} \right] \right. \nonumber \\
    &\quad \times \left. \left.  \prod_{i=k - \lambda+1}^{k-1} \exp\left\lbrace\frac{t}{\lambda} \left( \Vert \tilde r_i \Vert_2^2 - \ex\left[\Vert \tilde r_i \Vert_2^2 | \mathcal{F}_{i-1} \right] \right) \right\rbrace  \right\vert \siga \right] \\
    &\leq \mathbb{E} \Bigg[  \exp\left\lbrace \frac{t^2 \sigma^2 \bound^2}{2 \lambda^2} \right\rbrace\\
    &\quad \quad \times  \prod_{i=k - \lambda+1}^{k-1} \exp\left\lbrace\frac{t}{\lambda} \left( \Vert \tilde r_i \Vert_2^2 - \ex\left[\Vert \tilde r_i \Vert_2^2 | \mathcal{F}_{i-1} \right] \right) \right\rbrace  \Bigg\vert \siga \Bigg], \nonumber
    \end{eqnarray}
    \end{small}\
where we have made use of \cref{lemma:abs_diff} in the final line. Now, applying H\"{o}lder's inequality and the induction hypothesis,
\begin{small}
\begin{eqnarray}
    & \mathbb{E} \Bigg[ \exp\left\lbrace \frac{t^2 \sigma^2 \bound^2 }{2 \lambda^2} \right\rbrace\\
    &\quad \quad  \times \prod_{i=k - \lambda+1}^{k-1} \exp\left\lbrace\frac{t}{\lambda} \left( \Vert \tilde r_i \Vert_2^2 - \ex\left[\Vert \tilde r_i \Vert_2^2 | \mathcal{F}_{i-1} \right] \right) \right\rbrace  \Bigg\vert \siga \Bigg] \nonumber \\
    & \leq \mathbb{E}\left[ \left. \exp\left\lbrace \frac{t^2 \sigma^2 \bound^2}{2 \lambda} \right\rbrace \right\vert \siga \right]^\frac{1}{\lambda} \nonumber \\
    & \quad \times \mathbb{E} \left[ \left. \prod_{i=k - \lambda+1}^{k-1} \exp\left\lbrace\frac{t}{\lambda-1} \left( \Vert \tilde r_i \Vert_2^2 - \ex\left[\Vert \tilde r_i \Vert_2^2 | \mathcal{F}_{i-1} \right] \right) \right\rbrace  \right\vert \siga \right]^\frac{ \lambda - 1}{\lambda} \\
    & \leq \mathbb{E}\left[ \left. \exp\left\lbrace \frac{t^2 \sigma^2 \bound^2 }{2 \lambda} \right\rbrace \right\vert \siga \right]^\frac{1}{\lambda}\exp\left\lbrace \frac{t^2 \sigma^2 \bound^2 }{2(\lambda-1) } \sum_{j=1}^{\lambda-1} \frac{1}{j} \right\rbrace^\frac{ \lambda - 1}{\lambda}. \label{eqn-post-holder}
    \end{eqnarray}
\end{small}
    Since $\bound$ is measurable with respect to $\siga$, we have shown our desired inequality.
\end{proof}
From this distribution we can easily conclude the following corollary about the consistency of $\hat \rho_k^\lambda$.
\begin{corollary}\label{cor:cons-of-rho}
Suppose the setting of \cref{thm:convergence-of-moments} holds. Then, for any $\lambda \in \mathbb{N}$ and $\forall \delta > 0$, $\lim_{k \to \infty} \mathbb{P}( |\hat \rho_k^\lambda - \rho_k^\lambda | > \delta) = 0$. That is, $\hat \rho_k^\lambda$ is a consistent estimator of $\rho_k^\lambda$.
\end{corollary}
\begin{proof}
    We begin by noting that 
    \begin{small}
    \begin{eqnarray}
        \lim_{k \to \infty} \mathbb{P}( |\hat \rho_k^\lambda - \rho_k^\lambda | > \delta) &= \lim_{k \to \infty}\ex\left[\pr\left(|\hat \rho_k^\lambda - \rho_k^\lambda | > \delta|\siga\right)\right]\label{eq:tower} \\
        &\leq \lim_{k \to \infty}\ex \left[2e^{-\min\left\{\frac{\lambda \delta^2}{2\sigma^2\bound^2(1+\log(\lambda))},\frac{\delta}{2\omega\bound} \right\}}\right]\label{eq:bound_se}\\
        &= \ex \left[\lim_{k \to \infty} 2e^{-\min\left\{\frac{\lambda \delta^2}{2\sigma^2\bound^2(1+\log(\lambda))},\frac{\delta}{2\omega\bound} \right\}}\right]\label{eq:dom_con}\\
        &=0, \label{eq:conv_err}
    \end{eqnarray}
    \end{small}
where  (\ref{eq:tower}) comes from the tower property for conditional expectation;  (\ref{eq:bound_se}) comes from \cref{sub-Gaussian:estimator} and \cref{def:SE},  (\ref{eq:dom_con}) comes from $\bound$ being a monotonically decreasing sequence implying that the term inside the expectation is dominated by $2\exp(-\min\{\frac{\lambda \delta^2}{2\sigma^2M_0^2(1+\log(\lambda))},\frac{\delta}{2\omega M_0} \})$, an integrable function. This domination by an integrable function allows us to use the dominated convergence theorem to switch the limit and the integral.  (\ref{eq:conv_err}) then arises from \cref{thm:convergence-of-moments} and the definition of $\bound$.
\end{proof}
With the validity of our estimators established and the derivation of the distribution of $\hat \rho_k^\lambda$, we are now able to derive the uncertainty set and stopping criterion in the following two corollaries, whose proofs are straightforward.

\begin{corollary}\label{coll:credible-interval}
    An uncertainty set of level $1 - \alpha$ for $\rho_k^\lambda$
    is 
    \begin{equation} \label{eq:credible-interval-M}
    \hat  \rho_k^\lambda \pm \epsilon,
    \end{equation}
    where
    \begin{small}
    \begin{equation}
        \epsilon = \Bigg \{
        \begin{array}{ll}
            \sqrt{2 \log(2/\alpha)  \frac{\sigma^2 M^2_{k- \lambda}(1+\log(\lambda))}{\lambda}} & \text{if }\log(2 / \alpha) \leq \frac{\lambda \sigma^2 (1 + \log(\lambda))}{2\omega^2}\\
            \frac{2\log(2/\alpha)\omega\bound}{\lambda} & \text{if } \log(2 / \alpha) > \frac{\lambda \sigma^2 (1 + \log(\lambda))}{2\omega^2}. 
        \end{array} 
    \end{equation}
    \end{small} 
\end{corollary}

\begin{corollary}\label{corr:stopping}
    Given $\upsilon$, $\delta_I$, $\xi_I$, $\delta_{II}$, $\xi_{II}$ as defined in point \cref{subsection-solver-track-stop} and a sampling matrix satisfying \ref{assum:SE}, then:
    \begin{small}
    \begin{eqnarray} \label{eq:actual-typeI-control}
        &\bound \leq \min\left\{\frac{\lambda \upsilon^2(1 - \delta_I)^2}{2\log(1/\xi_I)\sigma^2 \bound (1+ \log(\lambda))},\frac{\lambda\upsilon(1-\delta_I)}{2\log(1/\xi_I) \omega }\right\}\\
        &\Rightarrow \mathbb{P}\left[ \hat \rho_{k}^\lambda > \upsilon, \rho_k^\lambda \leq \delta_I\upsilon \bigg{\vert} \mathcal{F}_{k-\lambda} \right] < \xi_I; \text{and}
        \end{eqnarray}
    \end{small}
    \begin{small}
        \begin{eqnarray}\label{eq:actual-typeII-control}
        &\bound \leq \min\left\{\frac{\lambda \upsilon^2(\delta_{II}-1)^2}{2\log(1/\xi_{II})\sigma^2 \bound (1+ \log(\lambda))},\frac{\lambda\upsilon(\delta_{II} - 1)}{2\log(1/\xi_{II}) \omega }\right\}\\
        &\Rightarrow \mathbb{P}\left[ \hat \rho_{k}^\lambda \leq  \upsilon, \rho_k > \delta_{II} \upsilon \bigg{\vert} \mathcal{F}_{k-\lambda} \right] < \xi_{II}.
        \end{eqnarray}
    \end{small}.
\end{corollary}

\subsection{Estimating the Uncertainty Set and Stopping Criterion}\label{subsec:est_stop}
\cref{coll:credible-interval,corr:stopping} provide a well-controlled uncertainty set and stopping criterion, yet require knowing $\bound^2$, which is usually not available. As stated before, \cref{coll:credible-interval,corr:stopping} can be operationalized by replacing $\bound^2$ with $\hat \iota_{k}^\lambda$. Of course, $\bound^2$ and $\hat \iota_k^\lambda$ must coincide in some sense in order for this estimation to be valid. Indeed, by \cref{thm:convergence-of-moments}, both $\bound^2$ and $\hat \iota_k^\lambda$ converge to zero as $k \to \infty$, which allows us to estimate $\bound^2$ with $\hat \iota_k^\lambda$ to generate consistent estimators. However, we could also estimate $\bound^2$ by $0$ to generate consistent estimators, but these would be uninformative during finite time. Therefore, we must establish that estimating $\bound^2$ by $\hat \iota_k^\lambda$ is also appropriate within some finite time. To do this we establish that the relative error between $\bound^2$ and $\hat \iota_k^\lambda$ is controlled by a constant (in probability). 

Accomplishing this task requires the following lemma, which gives a probability bound between $\hat \iota_k^\lambda$ and a specific intermediate quantity. 
\begin{lemma}\label{lem:rel-of-iota}
    Under the conditions of \cref{thm:convergence-of-moments}, $\forall \delta > 0$  
    \begin{small}
    \begin{eqnarray}
        &\pr\left(\left|\frac{\hat \iota_k^\lambda - \frac{\sum_{i = k - \lambda +1}^k (\ex[\|\tilde r_i \|_2^2 | \mathcal{F}_{i-1}])^2}{\lambda}}{\bound^2}\right| > \delta\bigg| \siga\right)\nonumber\\ 
         &\leq (1+\lambda)\exp\left( - \min\left\{\frac{\delta^2 \lambda}{2(2 + \sqrt{\delta \lambda})^2 \sigma^2(1+\log(\lambda))}, \frac{\delta \lambda}{2\omega (2+ \sqrt{\delta \lambda})}\right\}\right)
    \end{eqnarray}
    \end{small} 
\end{lemma}
\begin{proof}
    Using the definition of $\hat \iota_k^\lambda$ we have
    \begin{small}
    \begin{eqnarray}
        &\pr\left(\left|\frac{\hat \iota_k^\lambda - \frac{\sum_{i = k - \lambda +1}^k (\ex[\|\tilde r_i \|_2^2 | \mathcal{F}_{i-1}])^2}{\lambda}}{\bound^2}\right| > \delta\bigg| \siga\right)\nonumber\\
         &\quad \quad=\pr\left(\left|\sum_{i=k-\lambda + 1}^k\frac{ \|\tilde r_i\|_2^4 - (\ex[\|\tilde r_i \|_2^2 | \mathcal{F}_{i-1}])^2}{\lambda\bound^2} \right|> \delta \bigg| \siga\right) \nonumber \\
        &\quad \quad \leq \pr\left(\sum_{i=k-\lambda + 1}^k \left|\frac{ \|\tilde r_i\|_2^4 - (\ex[\|\tilde r_i \|_2^2 | \mathcal{F}_{i-1}])^2}{\lambda\bound^2} \right|> \delta \bigg| \siga\right) \nonumber \\
        & \quad \quad  \leq \pr\Bigg(\sum_{i=k-\lambda + 1}^k \left|\frac{ \|\tilde r_i\|_2^2 - \ex[\|\tilde r_i \|_2^2 | \mathcal{F}_{i-1}]}{\lambda\bound} \right|\nonumber\\ 
        &\quad \quad \quad \quad\quad \quad\times \left|\frac{ \|\tilde r_i\|_2^2 + \ex[\|\tilde r_i \|_2^2 | \mathcal{F}_{i-1}]}{\bound}\right|> \delta \bigg| \siga\Bigg)\label{eq:whole1}.
    \end{eqnarray}
    \end{small}
    Then by defining a variable, $G > 2$, to partition  (\ref{eq:whole1}) into disjoint sets and using the definition of measure,
    \begin{small}
    \begin{eqnarray}
        &\pr\left(\sum_{i=k-\lambda + 1}^k \left|\frac{ \|\tilde r_i\|_2^2 - \ex[\|\tilde r_i \|_2^2 | \mathcal{F}_{i-1}]}{\lambda\bound} \right| \left|\frac{ \|\tilde r_i\|_2^2 + \ex[\|\tilde r_i \|_2^2 | \mathcal{F}_{i-1}]}{\bound}\right|> \delta \bigg| \siga\right) \nonumber\\  
        &\quad =\pr\bigg(\sum_{i=k-\lambda + 1}^k \left|\frac{ \|\tilde r_i\|_2^2 - \ex[\|\tilde r_i \|_2^2 | \mathcal{F}_{i-1}]}{\lambda\bound} \right| \left|\frac{ \|\tilde r_i\|_2^2 + \ex[\|\tilde r_i \|_2^2 | \mathcal{F}_{i-1}]}{\bound}\right|> \delta , \nonumber\\ 
         &\quad \quad \quad \bigcap_{i=k-\lambda + 1}^k \left\{\left|\frac{ \|\tilde r_i\|_2^2 + \ex[\|\tilde r_i \|_2^2 | \mathcal{F}_{i-1}]}{\bound}\right|\leq G\right\} \bigg| \siga\bigg) \nonumber \\
         & \quad+\pr\bigg(\sum_{i=k-\lambda + 1}^k \left|\frac{ \|\tilde r_i\|_2^2 - \ex[\|\tilde r_i \|_2^2 | \mathcal{F}_{i-1}]}{\lambda\bound} \right| \left|\frac{ \|\tilde r_i\|_2^2 + \ex[\|\tilde r_i \|_2^2 | \mathcal{F}_{i-1}]}{\bound}\right|> \delta , \nonumber\\ 
         &\quad \quad \quad \bigcup_{i=k-\lambda + 1}^k \left\{\left|\frac{ \|\tilde r_i\|_2^2 + \ex[\|\tilde r_i \|_2^2 | \mathcal{F}_{i-1}]}{\bound}\right|> G\right\} \bigg| \siga\bigg)\nonumber\\
         & \leq\pr\bigg(\sum_{i=k-\lambda + 1}^k \left|\frac{ \|\tilde r_i\|_2^2 - \ex[\|\tilde r_i \|_2^2 | \mathcal{F}_{i-1}]}{\lambda\bound} \right| > \frac{\delta}{G} \bigg| \siga \bigg) \nonumber\\ 
         & \label{int-bound1} \quad \quad \quad + \pr\bigg(\bigcup_{i=k-\lambda + 1}^k \left\{\left|\frac{ \|\tilde r_i\|_2^2 + \ex[\|\tilde r_i \|_2^2 | \mathcal{F}_{i-1}]}{\bound}\right|> G\right\} \bigg| \siga\bigg).
    \end{eqnarray}
    \end{small}
    From here we will present the bounds for the left and right terms of  (\ref{int-bound1}) separately. For the left term of  (\ref{int-bound1}) we use  (\ref{eq:chernoff}) and when $t\leq \lambda / \omega$ have
    \begin{small}
    \begin{eqnarray}
      &\pr\bigg(&\sum_{i=k-\lambda + 1}^k \left|\frac{ \|\tilde r_i\|_2^2 - \ex[\|\tilde r_i \|_2^2 | \mathcal{F}_{i-1}]}{\lambda \bound} \right| > \frac{\delta}{G}\bigg| \siga\bigg)  \label{rchernoff1}\\
        &&\leq  \exp \left(\frac{t^2 \sigma^2 (1+ \log(\lambda))}{2\lambda} - \frac{\delta t}{G}\right)\\
        &&\label{rmint1}\leq \exp\left( - \min\left\{\frac{\delta^2 \lambda}{2G^2 \sigma^2(1+\log(\lambda))}, \frac{\lambda\delta}{2\omega G}\right\}\right), 
    \end{eqnarray}
    \end{small}
    where  (\ref{rmint1}) comes from minimizing  (\ref{rchernoff1}) in terms of $t$ under the constraint $t\leq \lambda / \omega$. We next address the right side of  (\ref{int-bound1}) for which we have
    \begin{small}
    \begin{eqnarray}
        \pr&\bigg(\bigcup_{i=k-\lambda + 1}^k \left\{\left|\frac{ \|\tilde r_i\|_2^2 + \ex[\|\tilde r_i \|_2^2 | \mathcal{F}_{i-1}]}{\bound}\right|> G\right\} \bigg| \siga\bigg)  \nonumber \\
         & =\pr\bigg(\bigcup_{i=k-\lambda + 1}^k \left\{\left|\frac{ \|\tilde r_i\|_2^2 - \ex[\|\tilde r_i \|_2^2 | \mathcal{F}_{i-1}] + 2 \ex[\|\tilde r_i \|_2^2 | \mathcal{F}_{i-1}]}{\bound}\right|> G \right\}\bigg| \siga\bigg)\nonumber\\
        &\label{bounddef1} \leq \pr\bigg(\bigcup_{i=k-\lambda + 1}^k \left\{\left|\frac{ \|\tilde r_i\|_2^2 - \ex[\|\tilde r_i \|_2^2 | \mathcal{F}_{i-1}]}{\bound}\right| + 2> G\right\} \bigg| \siga\bigg)  \\
        &  \leq \sum_{i=k-\lambda + 1}^k\pr\bigg(\left|\frac{ \|\tilde r_i\|_2^2 - \ex[\|\tilde r_i \|_2^2 | \mathcal{F}_{i-1}]}{\bound}\right|> G - 2 \bigg| \siga\bigg)\\
        &\label{lchernoff1}   \leq \lambda \exp \left(\frac{t^2\sigma^2 }{2} - t\left(G - 2\right)\right)\\
        & \label{mincher1}  \leq \lambda \exp\left( -\min\left\{\frac{\left(G - 2\right)^2}{2\sigma^2}, \frac{G-2}{2\omega}\right\} \right),
    \end{eqnarray}
\end{small}
    where  (\ref{bounddef1}) comes from  (\ref{rtoM}),
     (\ref{lchernoff1}) comes from the  (\ref{eq:chernoff}), and  (\ref{mincher1}) comes from minimizing  (\ref{lchernoff1}) in terms of $t$ under the constraint $t\leq 1 / \omega$. Putting both parts together gives us,
    \begin{small}
    \begin{eqnarray*}
        &\pr\left(\left|\frac{\hat \iota_k^\lambda - \iota_k^\lambda}{\bound^2}\right| > \delta\bigg| \siga\right) \\
        &\quad \quad\leq \inf_{G > 2} \exp\left( - \min\left\{\frac{\delta^2 \lambda}{2G^2 \sigma^2(1+\log(\lambda))}, \frac{\delta\lambda}{2\omega G}\right\}\right) + \\ &\quad \quad \quad\quad\lambda \exp\left( -\min\left\{\frac{\left(G - 2\right)^2}{2\sigma^2}, \frac{G-2}{2\omega}\right\} \right). \label{final-rate-iota1}
    \end{eqnarray*}       
    \end{small}
    We can then observe that when $G \geq 2 + \sqrt{\delta \lambda}$ it is the case that
    \begin{eqnarray*}
        \exp&\left( - \min\left\{\frac{\delta^2 \lambda}{2G^2 \sigma^2(1+\log(\lambda))}, \frac{\delta \lambda}{2\omega G}\right\}\right)\\
         &\geq \exp\left( -\min\left\{\frac{\left(G - 2\right)^2}{2\sigma^2}, \frac{G-2}{2\omega}\right\} \right). 
    \end{eqnarray*}
    Thus, if we let $Y = 2 + \sqrt{\delta \lambda}$ we can upper bound the right-hand side of  (\ref{final-rate-iota1}) in the following manner,
    \begin{small}
    \begin{eqnarray*}
        &\inf_{G > 2} \exp\left( - \min\left\{\frac{\delta^2 \lambda}{2G^2 \sigma^2(1+\log(\lambda))}, \frac{\delta \lambda}{2\omega G}\right\}\right)\\&\quad + \lambda\exp\left( -\min\left\{\frac{\left(G - 2\right)^2}{2\sigma^2}, \frac{G-2}{2\omega}\right\} \right)\\
        &\leq \inf_{G > Y} (1+\lambda)\exp\left( - \min\left\{\frac{\delta^2 \lambda}{2G^2 \sigma^2(1+\log(\lambda))}, \frac{\delta \lambda}{2\omega G}\right\}\right)\\
        & \leq (1+\lambda)\exp\left( - \min\left\{\frac{\delta^2 \lambda}{2(2 + \sqrt{\delta \lambda})^2 \sigma^2(1+\log(\lambda))}, \frac{\delta \lambda}{2\omega (2+ \sqrt{\delta \lambda})}\right\}\right).
    \end{eqnarray*}
    \end{small}
\end{proof}
We now present our bound on the relative error between $\hat \iota_k^\lambda$ and $\bound^2$.

\begin{theorem} \label{lemm:relative-error-M-iota}
    Under the conditions of \cref{thm:convergence-of-moments}, for $\delta > 0$, $\bound^2$ as described in \cref{sub-Gaussian:estimator},
    \begin{small} 
    \begin{eqnarray}
    &\mathbb{P}\left( \left| \frac{\bound^2 - \hat \iota_k^\lambda}{\bound^2}\right| > 1 + \delta, \bound^2 \neq 0  \bigg| \siga \right)\\ 
    &\leq (1+\lambda)\exp\left( - \min\left\{\frac{\delta^2 \lambda}{2(2 + \sqrt{\delta \lambda})^2 \sigma^2(1+\log(\lambda))}, \frac{\delta \lambda}{2\omega (2+ \sqrt{\delta \lambda})}\right\}\right).
    \end{eqnarray}
    \end{small}
\end{theorem}
\begin{proof}
First,
\begin{small}
\begin{eqnarray}
    \left| \frac{\bound^2 - \hat \iota_k^\lambda}{\bound^2}\right| &
    \leq \left| \frac{\bound^2 - \frac{\sum_{i = k - \lambda +1}^k (\ex[\|\tilde r_i \|_2^2 | \mathcal{F}_{i-1}])^2}{\lambda}}{\bound^2}\right| \\
    &\quad +  \left| \frac{\frac{\sum_{i = k - \lambda +1}^k (\ex[\|\tilde r_i \|_2^2 | \mathcal{F}_{i-1}])^2}{\lambda} - \hat \iota_k^\lambda}{ \bound^2} \right| \nonumber\\
    &\leq 1 
    + \left| \frac{\frac{\sum_{i = k - \lambda +1}^k (\ex[\|\tilde r_i \|_2^2 | \mathcal{F}_{i-1}])^2}{\lambda} - \hat \iota_k^\lambda}{\bound^2} \right|.\label{fact-bound1} 
\end{eqnarray}
\end{small}
We now apply the bound in \cref{lem:rel-of-iota} to conclude.
\end{proof}

Owing to \cref{lemm:relative-error-M-iota}, the relative error between $\hat \iota_k^\lambda$ and $\bound^2$ is reasonably well controlled for practical purposes. As a result,  we can use $\hat \iota_k^\lambda$ as a plug-in estimator for $\bound^2$ for the uncertainty set, \eref{eq:credible-interval-M}, to produce the estimated uncertainty set suggested in Line \ref{Credible-Interval} of \cref{algo:Adaptive-window}. Additionally, we can use  $\hat \iota_k^\lambda$ as a plug-in estimator for $\bound^2$ for the stopping condition controls in  (\ref{eq:actual-typeI-control}) and  (\ref{eq:actual-typeII-control}) to produce the estimated stopping criterion in Line \ref{stopping-condition} of \cref{algo:Adaptive-window}.

\section{Experimental results}\label{sec:exp}
We conduct two experiments.  Our first experiment is solving a collocation problem reformulated as a streaming problem (see \ref{ex:colloc}). We close this section with a timing comparison between \cref{algo:Adaptive-window} and periodic residual calculation methods on a 125,000 by 125,000 system.

\subsection{Collocation Problem}\label{subsec:boundary}
To demonstrate that \cref{algo:Adaptive-window} works on large-scale streaming problems, we turn to the collocation problem laid out in \ref{ex:colloc}. In this problem, we verify that the uncertainty set for $\hat{\rho}_k^\lambda$ contains $\rho_k^\lambda$ at appropriate rates and that the stopping criterion correctly controls stopping errors. To accomplish this task we perform the following experiment.
\begin{experiment}\label{Exp:Coll}
    Given the gap between fixed grid points $\epsilon = 1/99$, the number of sample rows $p=20$, the  moving average width $\lambda_1\in \{100,300\}$, the uncertainty set parameter $\alpha$, and the stopping criterion parameters $(\upsilon,\delta_I,\delta_{II}, \xi_I, \xi_{II}) = (400, .9, 1.1, .01,.01)$ we: 
\begin{enumerate}
    \item Generate an equally spaced grid on the unit cube with the gap between grid points being $\epsilon$. This equates to there being $(1+1/\epsilon)^3$ fixed points on the cube. 
    \item Generate $(\tilde A_k,\tilde b_k)$ by generating 20 random coordinates as explained in \ref{ex:colloc}. 
    \item Run \cref{algo:Adaptive-window} using $(\tilde A_k,\tilde b_k)$'s and save at each iteration the $x_k$, $\hat \rho_k^\lambda$, $\hat \iota_k^\lambda$, $\|r_k\|_2^2$, and the width of moving average $\lambda$.
    \item At each $x_k$ approximate $\ex[\|\tilde r_k\|_2^2 | \mathcal{F}_{k-1}]$ using a Monte Carlo simulation with 100 samples.
    \item Approximate $\rho_k^\lambda$ at the $k^{\text{th}}$ iteration by taking a moving average of the $\hat \ex[\|\tilde r_k\|_2^2| \mathcal{F}_{k-1}]$s with the same width as $\hat \rho_k^\lambda$ at iteration $k$. 
\end{enumerate}
\end{experiment}

Practically, our value of $\sigma^2$ in \eref{colloc:exp} (i.e., $V/2$) is roughly 9 trillion in this experiment, which is far too conservative to be useful.\footnote{Based on different choices of $\epsilon$ it does not appear the $V/2$ ever really accurately reflects the variance.} Thus, instead of using $V/2$ for $\sigma^2$, we estimate it by computing
\begin{small} 
\begin{equation}
    \text{Variance}\left(\frac{|\hat \ex[\|\tilde r_k\|_2^2| \mathcal{F}_{k-1}] - \|r_k\|_2^2|}{\hat \ex[\|\tilde r_k\|_2^2| \mathcal{F}_{k-1}]}, k = 1, \cdots, 125 \right),
\end{equation}\end{small}
with $\|r_k\|_2^2$ coming from step (iii) and $\hat \ex[\|\tilde r_k\|_2^2| \mathcal{F}_{k-1}]$ coming from step (v) of \cref{Exp:Coll}. This estimator is more optimistic than $V/2$, yet seems to be appropriate. It should be noted that the choice of using the first 125 iterations for the variance estimation is somewhat arbitrary and many other choices of the number iterations could be made to produce approximately the same variance estimate as can be seen in \tref{tab:failure-k}.

The results of this experiment can be observed in \fref{fig:big} for two different choices of $\lambda_1$ one being $\lambda_1 = 100$ and one being $\lambda_1=300$. As should be clear from \fref{fig:big} the $95\%$ uncertainty sets correctly cover $\rho_k^\lambda$ at $99.6\%$ of the iterations for $\lambda_1 = 300$ and $99.4\%$ for $\lambda_1 = 100$, which is still conservative despite us not using the larger estimate of $\sigma^2$ from \eref{colloc:exp}. This rate coverage failure also does not change much when the width is set to be some other value as seen in \tref{tab:failure}. Additionally, when the stopping criterion is satisfied there are no instances where the two types of stopping errors occur. This indicates that even at large scales \cref{algo:Adaptive-window} still performs well.

\begin{figure}
    \centering
   \begin{tikzpicture}
    \begin{axis}[
        ymode = log,
        width = .45 \textwidth,
        height = .20 \textheight,
        every axis plot/.append style={ultra thick},
        xlabel = {Iteration},
        ylabel = {Norm Squared},
        title = {Uncertainty set results for $\rho_k^\lambda$ in Collocation Problem},
        legend style = {font = \tiny},
        label style={font=\tiny},
        title style={font=\tiny},
        ylabel style = {yshift = -.25cm},
        xlabel style = {yshift = .2cm},  
        tick label style={font=\tiny}]
        \addplot[black] table[x = it, y = Upper, col sep=comma, mark = none]{\tpath/big.csv}; 
        \addplot[black] table [x = it, y = Lower, col sep=comma, mark = none]{\tpath/big.csv};
        \addplot[red!60!black, mark = none] table [x = it, y = Mov_Obs_Res, col sep=comma]{\tpath/big.csv};
        \addplot[green!60!black, mark = none] table [x = it, y = Mov_True_Res, col sep=comma]{\tpath/big.csv};
        
        \legend{,,$\hat \rho_k^\lambda$, $\rho_k^\lambda$}
    \end{axis}
\end{tikzpicture}
%
   \begin{tikzpicture}
    \begin{axis}[
        ymode = log,
        width = .45 \textwidth,
        height = .20 \textheight,
        every axis plot/.append style={ultra thick},
        xlabel = {Iteration},
        ylabel = {Norm Squared},
        title = {Uncertainty set results for $\rho_k^\lambda$ in Collocation Problem},
        legend style = {font = \tiny},
        label style={font=\tiny},
        title style={font=\tiny},
        ylabel style = {yshift = -.25cm},
        xlabel style = {yshift = .2cm},  
        tick label style={font=\tiny}]
        \addplot[black] table[x = it, y = Upper, col sep=comma, mark = none]{\tpath/big100.csv}; 
        \addplot[black] table [x = it, y = Lower, col sep=comma, mark = none]{\tpath/big100.csv};
        \addplot[red!60!black, mark = none] table [x = it, y = Mov_Obs_Res, col sep=comma]{\tpath/big100.csv};
        \addplot[green!60!black, mark = none] table [x = it, y = Mov_True_Res, col sep=comma]{\tpath/big100.csv};
        
        \legend{,,$\hat \rho_k^\lambda$, $\rho_k^\lambda$}
    \end{axis}
\end{tikzpicture}
%
  \caption{In the left plot, the black lines represent the upper and lower bounds of the uncertainty set, the green line represents $\rho_k^\lambda$, the red line represents $\hat \rho_k^\lambda$ for $\lambda = 300$ and the right plot represents the same values for $\lambda = 100$.} \label{fig:big}
\end{figure}

\begin{table}
    \caption{Interval failure rates at different moving average window widths}
\centering
 \begin{tabular}{|c|c|c|c|c|c|}
    \hline 
    Width & 10 & 50 & 100 & 200 &300\\
    \hline
    Failure Rate& 0.0148& 0.0148&0.006&0.004&0.004\\
    \hline
 \end{tabular}  \label{tab:failure} 
\end{table}

\begin{table}
    \caption{Interval failure rates with different $\sigma^2$ estimates and moving average window widths}
    \centering
    \begin{tabular}{|l|l|l|l|l|l|l|}
      \hline
      \multirow{2}{*}{Dataset} &
        \multicolumn{2}{c|}{$k=5$} &
        \multicolumn{2}{c|}{$k=125$} &
        \multicolumn{2}{c|}{$k=474$} \\
        \cline{2-7}
      & $\lambda = 10$ & $\lambda = 300$ & $\lambda = 10$ & $\lambda = 300$ &$\lambda = 10$ & $\lambda = 300$ \\
      \hline
      $\sigma^2$ & 0.172 & 0.172 & 0.079 & 0.079 & 0.088 & 0.088 \\
      \hline
      Failure Rate & 0.0148 & 0.004& 0.0148 & 0.004 & 0.0148 & 0.004\\
      \hline
    \end{tabular}\label{tab:failure-k}
  \end{table}
  
\subsection{Timing comparison with naive methods}\label{subsec:time-comparison}
We now demonstrate the computational benefits of our method over periodic calculation of the full residual on a  125,000 by 125,000 system. We performed this experiment using a single thread of a Xeon E5-2680 v3 @ 2.50GHz with 32 GB of RAM. The results of the experiment can be found in \tref{table-solvetimes}.

\begin{table}[!h]
    \caption{Iterations achieved in four-day run on a collocation problem with a $51 \times 51 \times 51$ grid collocation points}
   \label{table-solvetimes}
    \centering
    \begin{tabular}{llll}
        \toprule
        \textbf{Method} & \textbf{Iterations} &  \textbf{Time per iteration}\\
        \midrule
          \Cref{algo:Adaptive-window} & $7819$ & $44.2$ seconds\\
          Full Residual every GBRK update & $16$ & $22,109$ seconds\\ 
          Full Residual every 1000 GBRK updates & $4162$ & $83.03$ seconds\\
        \bottomrule
    \end{tabular} 
  \end{table}
  
When comparing the performance of \cref{algo:Adaptive-window} to other methods of tracking, we can see substantial benefits of \cref{algo:Adaptive-window} compared to alternatives. We do this by calculating the number of iterations completed in a four day window by three methods: \cref{algo:Adaptive-window}, calculating the full residual after each GBRK update, and calculating the full residual after every 1000 GBRK updates.

\section{Conclusion}\label{sec:Con}
To address the issue of effectively tracking and stopping the progress of a streaming solver, we have presented a computationally efficient estimator and uncertainty set for a moving average of the residuals. We then rigorously demonstrated the effectiveness of this estimator using only assumptions about the streams having a set of consistent solutions and the norm of the residuals of the streams being sub-Exponential.  We additionally show that the assumptions used to make these conclusions are relatively weak and apply to a large class of common problems. Moreover, we verified our methodology by successfully applying it to a large-scale collocation problem and demonstrated its computational benefits over alternative methods.

\ack{This material is based upon work supported by the National Science Foundation (NSF) under grant no. 2309445.}

\section*{References}
\bibliographystyle{unsrt}
\bibliography{../bib/linear_res}
\appendix
\section{Areas of Applications}\label{sec:app}
With the theory for \cref{algo:Adaptive-window} established, we now move on to presenting common areas of application that satisfy \cref{assum:consis,assum:SE} and for which it is sensible apply \cref{algo:Adaptive-window}. The situations that will be discussed are Johnson-Lindenstrauss matrix sketches, Block Randomized Kaczmarz, and the collocation problem arising from boundary element analysis.

\subsection{Johnson-Lindenstrauss Matrix Sketching}\label{ex:diff_priv}
With the pervasiveness of large scale linear system in areas such as Gaussian Process modeling \cite{Ramussen}, optimization \cite{dembo}, and network analysis \cite{zhan_identification_2017}, random algorithms have been shown to accelerate the speed of achieving good approximate solutions. One popular form of random transformation is the use of a random matrix that satisfies the Johnson-Lindenstrauss property \cite{ACHLIOPTAS2003671,Ailon2009TheFJ,JL}, whose definition follows.

\begin{definition}\label{def:JL}
    A matrix $S \in \mathbb{R}^{m \times p}$ satisfies the Johnson-Lindenstrauss property if there exists constants $C,\omega > 0$ s.t. for all $\delta \geq 0$ and for any $x \in \mathbb{R}^m$,
    \begin{eqnarray}
         \mathbb{P}\left(|\|S^\top x\|_2^2 - \|x\|_2^2| > \delta \|x\|_2^2\right) \quad < 2e^{-\min\left\{(Cp\delta^2)/2,\delta/(2\omega) \right\}}.
    \end{eqnarray}
\end{definition}

\begin{remark}
Matrices $S$ satisfying \cref{def:JL} are known as sketching matrices and have broad applications in numerical linear algebra \cite{drineas2016randnla,martinsson2020randomized}. Examples of such matrices include the Gaussian matrix \cite{indyk}, the Achlioptas sparse sampling matrix \cite{ACHLIOPTAS2003671}, and the Fast Johnson-Lindenstrauss transform \cite{Ailon2009TheFJ}.  These methods will be the focus of experiment section, and thus we have included the values $C$ and $\omega$ for when $\delta = 1$ in \tref{table-JL-constant-values}.
\end{remark}
By applying a matrix satisfying \cref{def:JL} to a consistent linear system with coefficient matrix $A \in \mathbb{R}^{m \times n}$ and constant vector $b \in \mathbb{R}^m$, and we denote $\tilde A = S^\top A$ and $\tilde b = S^\top b$. Then, we readily see that $(\tilde A, \tilde b)$ satisfy \cref{assum:SE,assum:consis} \cite{Gupta_2017}. This type of process, in addition to allowing faster approximations of large consistent systems, has the potential for use in large network systems where one may compute an Adjacency-based, or Graph-Laplacian-based, measure such as Katz centrality \cite{zhan_identification_2017}, and still maintain the privacy of the individual according to a particular differential privacy standard \cite{updadhyayDP,blockiDP}. In this case, one could use realizations a of Johnson-Lindenstrauss sketch and generate $(\tilde{A}_k,\tilde{b}_k)$, which can then be used in the Generalized Block Randomized Kaczmarz (GBRK) framework to compute the desired metric without violating the privacy of the individuals.
\begin{table}[h]
    \centering
    \caption{Values of $C$ and $\omega$ in \cref{def:JL} for common sampling methods.}
    \label{table-JL-constant-values}
    \begin{tabular}{lccccc}
          \toprule
           &$C$ & $\omega$\\
           \midrule
          Gaussian Matrix \cite{Dasgupta2003AnEP}  & 1.1 & .47\\
           Achlioptas \cite{ACHLIOPTAS2003671}  & 1.16 & .46\\
           FJLT \cite{Ailon2009TheFJ} & 0.83 & .7\\ 
           \bottomrule
       \end{tabular} 
\end{table}
\begin{table}[!hb]
    \caption{Conservative choices of the contraction parameter, $\eta$, by sketching method.\protect\footnotemark}
    \label{table-eta}
    \centering
    \begin{tabular}{llll}
        \toprule
        \textbf{Sampling Method} & Achlioptas \cite{ACHLIOPTAS2003671} & FJLT \cite{Ailon2009TheFJ}  & Gaussian \cite{Dasgupta2003AnEP}\\
        \midrule
          $\eta$ & $26$ & $188$ & $26$\\
        \bottomrule
    \end{tabular} 
  \end{table}\footnotetext{By conservative, we mean that these are choices of $\eta$ that are the largest integer such that the coverage failure rate of the interval is as large as possible, but still less than the designed rate.} 
\subsection{Random Subsets (Randomized Block Kaczmarz)}\label{ex:rand_sub}
One of the simplest versions of the streaming problem is where the observations are generated from random subsets of some larger system. Specifically, a sequence of observations, $\lbrace (\tilde A_k,\tilde b_k): k \in \mathbb{N} \rbrace$, that are independent, identically distributed subsets of the equations of a larger consistent system, whose coefficient matrix is given by $A \in \mathbb{R}^{m \times n}$ and whose constant vector is $b \in \mathbb{R}^m$. This is often specifically referred to as the Randomized Block Kaczmarz method \cite{Necora}. This type of method is often used for solving tomography problems under the name Algebraic Reconstruction Technique \cite{birk2014gpu,kakCT}.

This model satisfies \cref{assum:SE}, even with the allowances that each row is weighted by its own respective constant and each block, $\mathcal{J}$, has its own probability of selection given by $P_\mathcal{J}$ (as this is typically how these block Kaczmarz methods are implemented \cite{Necora}). This is stated formally now.

\begin{proposition}\label{prop:submat}
    Let $A \in \mathbb{R}^{m \times n}$ and $b \in \mathbb{R}^m$ such that $\exists x^* \in \mathbb{R}^n$ where $Ax^* = b$.
        Let $\mathcal{J}$ be a random subset of $\lbrace 1,\ldots,m \rbrace$ and let $P_i > 0$ be the probability $ i \in \mathcal{J}$. Let $w_i > 0$ for $i=1,\ldots,m$. Define $S_{\mathcal{J}} \in \mathbb{R}^{m \times p}$ to be the matrix whose columns are scalings of standard basis elements indexed by $\mathcal{J}$, where the scaling for the column corresponding to index $i \in \mathcal{J}$ is $w_i$. Let $(\tilde A_{\mathcal{J}} , \tilde b_{\mathcal{J}}) = (S_{\mathcal{J}}^\top A,S_{\mathcal{J}}^\top b)$. Then, for any $x \in \mathbb{R}^n$ such that $Ax \neq b$, 
        \begin{equation} \label{row_subset:exp}
            \frac{\|\tilde{A}_{\mathcal{J}} x - \tilde b_{\mathcal{J}}\|^2_2 - \ex [\|\tilde A_{\mathcal{J}} x - \tilde b_{\mathcal{J}}\|^2_2]}{ \ex [\|\tilde A_{\mathcal{J}} x - \tilde b_{\mathcal{J}}\|^2_2]} \sim \text{SE}\left(\frac{\max_i w_i^2}{2(\min_{i} P_{i} w_i^2)},{0}\right).
        \end{equation}
    \end{proposition}
\begin{proof}
        By \cref{lemma:boundedSE}, all bounded random variables are sub-Exponential with $\omega = 0$. To show the block residuals are bounded in this instance, we find a lower bound on $\ex [ \Vert \tilde A_{\mathcal{J}} x - \tilde b_{\mathcal{J}} \Vert_2^2]$ and an upper bound on $\|\tilde A_{\mathcal{J}} x - \tilde b_{\mathcal{J}}\|_2^2$. For the lower bound, let $a_i \in \mathbb{R}^n$ denote the $i^{\text{th}}$ row of $A$, and let $\mathbb{P}(i \in \mathcal{J}) = P_{i}$. This allows us to write 
        \begin{eqnarray}
            \ex [ \Vert \tilde A_{\mathcal{J}} x - \tilde b_{\mathcal{J}} \Vert_2^2] &= \sum_{i = 1}^m P_{i} w_i^2 [a_i^\top (x - x^*)]^2 \\ 
            &\geq (\min_i P_i w_i^2 ) \sum_{i = 1}^m [a_i^\top (x - x^*)]^2 
        \end{eqnarray}

For the upper bound on $\|\tilde A_{\mathcal{J}} x - \tilde b_{\mathcal{J}}\|_2^2$ we note that 
\begin{equation}
    \sum_{i\in \mathcal{J}} w_i^2 [a_i^{\top}(x - x^*)]^2 \leq (\max_i w_i^2) \sum_{i = 1}^m [a_i^{\top}(x - x^*)]^2.
\end{equation}
       
       Combing these lower and upper bounds, we conclude that for any $x$ such that $Ax \neq b$,
\begin{equation}
-1 \leq \frac{\|\tilde{A}_{\mathcal{J}} x - \tilde b_{\mathcal{J}}\|^2_2- \ex[\|\tilde A_\mathcal{J}x - \tilde b_\mathcal{J}\|_2^2]}{\ex [\|\tilde A_{\mathcal{J}} x - \tilde b_{\mathcal{J}}\|^2_2]} 
\leq \frac{\max_i w_i^2}{\min_i P_i w_i^2} - 1.
\end{equation}       
       
Giving us that the relative error is bounded. We then obtain the result by applying \cref{lemma:boundedSE}.
       \end{proof}
    
In the case where the rows are selected with equal probability and equal weights, this bound is equivalent to the variance bound being $m/p$. As this grows with the dimension of the matrix, concerns may arise relating to the tightness of this bound. While a tighter bound is not possible for general matrices, we can tighten this bound depending on the block condition numbers in the following proposition.

\begin{proposition}\label{prop:tighterUB}
    Let the conditions in \cref{prop:submat} hold. Further, if we let $s_{\max_p},s_{\min_p}$ be the largest and smallest singular value of all row blocks of a matrix $A$ containing $p$ rows. Then, for any $x \in \mathbb{R}^n$ such that $Ax \neq b$, 
            \begin{equation} \label{prop:tigheterUB:exp}
                \frac{\|\tilde{A}_{\mathcal{J}} x - \tilde b_{\mathcal{J}}\|^2_2 - \ex [\|\tilde A_{\mathcal{J}} x - \tilde b_{\mathcal{J}}\|^2_2]}{ \ex [\|\tilde A_{\mathcal{J}} x - \tilde b_{\mathcal{J}}\|^2_2]} \sim \text{SE}\left(\frac{s^2_{\max_{p}}}{2s^2_{\min_{p}}},{0}\right).
            \end{equation}
\end{proposition}
\begin{proof}
    We note that
   \begin{eqnarray}
    \frac{\|\tilde{A}_{\mathcal{J}} x - \tilde b_{\mathcal{J}}\|^2_2}{\ex[\|\tilde{A}_{\mathcal{J}} x - \tilde b_{\mathcal{J}}\|^2_2]} &\leq \frac{s^2_{\max_{p}}}{s^2_{\min_{p}}}.
   \end{eqnarray}
   This allows us to conclude that
   \begin{equation}
    -1 \leq \frac{\|\tilde{A}_{\mathcal{J}} x - \tilde b_{\mathcal{J}}\|^2_2- \ex[\|\tilde A_\mathcal{J}x - \tilde b_\mathcal{J}\|_2^2]}{\ex [\|\tilde A_{\mathcal{J}} x - \tilde b_{\mathcal{J}}\|^2_2]} 
    \leq \frac{s^2_{\max_{p}}}{s^2_{\min_{p}}}  - 1.
    \end{equation} 
    The result then follow from \cref{lemma:boundedSE}.      
\end{proof}


\subsection{Collocation Problem with Random Coordinates} \label{ex:colloc}
Our final example of a problem that satisfies \cref{assum:SE} is the collocation problem in boundary element analysis. We will describe this problem in detail for a specific setting that will be used in our largest numerical experiment (see \cref{subsec:boundary}); however, the results from this subsection will hold for more general problems as well.

 In this problem, we look to approximate the unknown solution, $u(t)$, with $t = (t^{(1)},t^{(2)},t^{(3)}) \in [0,1]^3$, that satisfies 
\begin{small} 
 \begin{eqnarray}\label{eq:interior}
    \Delta(u(t)) = \frac{-7 \pi^2}{2} \sin(\pi t^{(1)})\sin\left(\frac{\pi t^{(2)}}{2}\right)\sin\left(\frac{3 \pi t^{(3)}}{2}\right), \text{ } t \in [0,1]^3,\\
    u(t) = \sin(\pi t^{(1)}) \sin(\pi t^{(2)}/2) \sin(3\pi t^{(3)}/2), \text{ } t \in \partial [0,1]^3,
\end{eqnarray}
\end{small}
where $t \in \partial [0,1]^3$ represents the $t$ being on boundary of the cube; and $\Delta$ is the Laplace operator. 

To make this approximation, we wish to find the coefficients, $x_j$, of the linear combination of Quadric Radial Basis functions, $\phi(t,\chi_j) = \sqrt{\Vert t - \chi_j\Vert_2^2 + 1}$, and their Laplacians, evaluated at a fixed and finite set of control points $\lbrace \chi_j \rbrace$ equally spaced throughout the unit cube.  

Since both the Laplacian and the boundary condition are linear operators, we can represent \ref{eq:interior} as a streaming linear system, whose streams consist of coordinates $t_i \in [0,1]^3$ selected uniformly at random from either the interior (with probability 2/3), the faces (with probability 1/6), or the edges (with probability 1/6). 

With the sampling points $t_i$ and control points $\chi_j$, we form a linear system defined by
\begin{small}
\begin{eqnarray}
    \tilde A^{(i,j)} &= \Bigg\{
    \begin{array}{ll}
        \phi_j(t_i,\chi_j) &\text{ if } t_i \in \partial [0,1]^3 \\
         \Delta  \phi_j(t_i,\chi_j) &\text{ if } t_i \in (0,1)^3 
    \end{array},\label{eq:collocA} \\
    \tilde b^{(i)} &= \Bigg\{\begin{array}{ll}
        \sin(\pi t^{(1)}_i) \sin(\pi t^{(2)}_i/2) \sin(3\pi t^{(3)}_i/2) &\text{ if } t_i \in \partial [0,1]^3\\
        \frac{-7 \pi^2}{2} \sin(\pi t^{(1)}_i)\sin\left(\frac{\pi t^{(2)}_i}{2}\right)\sin\left(\frac{3 \pi t^{(3)}_i}{2}\right) &\text{ if } t_i \in (0,1)^3.\label{eq:collocb}
    \end{array} 
\end{eqnarray}
\end{small}

When the set of sampling and control points are the same, the resulting linear system is consistent as $A$ is nonsingular \cite[Theorem 2.2]{buhmann2003radial}. Below, however, these two sets are allowed to differ and, in our experiments, the set of sampling points are drawn at random. This random set of sampling points brings into question whether \cref{assum:consis} is true. While this is often ignored in practice, this concern can relieved by two observations. First, if a substantial number of control points are used, then the approximation discrepancy between the target function, \ref{eq:interior}, and the linear combination of basis functions can be made sufficiently small. Second, if there is such a well-controlled discrepancy, then the procedure will generate a solution whose error is bounded by a scaling of this discrepancy \cite[Theorem 2.1]{needell2010randomized}. Hence, in our experiments below, we use a large number of control points to alleviate concerns about \cref{assum:consis}. Fortunately, we can show that this problem satisfies \cref{assum:SE}.

\begin{proposition}\label{prop:colloc}
    Let $\lbrace t_1,\ldots,t_p \rbrace$ be chosen independently, from $[0,1]^3$ as described above. Let the entries of $\tilde A_k$ be defined according to \ref{eq:collocA} and those of $\tilde b_k$ be defined according to \ref{eq:collocb}. Then, if we let $V = {9N}/{\sigma_{\min}( \ex[ \tilde A_k^\top \tilde A_k])}$, where $\sigma_{\min}$ is the smallest non-zero singular value and $N$ is the number of columns in $\tilde A_k$, 

    \begin{equation} \label{colloc:exp}
        \frac{\|\tilde{A}_kx - \tilde b_k\|^2_2 - \ex [\|\tilde A_k x - \tilde b_k\|^2_2]}{ \ex [\|\tilde A_k x - \tilde b_k\|^2_2]} \sim \text{SE}\left(\frac{V}{2},{0}\right). 
    \end{equation}
\end{proposition}
\begin{proof}

    Let $x^* \in \mathbb{R}^n$ denote any vector such that $\mathbb{P}(\tilde A_k x^* = \tilde b_k) = 1$. For any $x \in \mathbb{R}^n$ such that $\mathbb{P}( \tilde A_k x = \tilde b_k) < 1$, we can decompose $x - x^*$ into $u \in \mathrm{row}( \ex [ \tilde A_k^\top \tilde A_k ])$ and $v \in \mathrm{null}( \ex [\tilde A_k^\top \tilde A_k ])$. By construction, $u \neq 0$. Moreover, $v \in \mathrm{null}(\tilde A_k^\top \tilde A_k)$ with probability one (otherwise we would have a contradiction with $v \in \mathrm{null}( \ex [\tilde A_k^\top \tilde A_k ])$).
    
    Therefore,
    \begin{equation}
    0 \leq \frac{ \Vert \tilde A_k x - \tilde b_k \Vert_2^2 }{ \ex[ \Vert \tilde A_k x - \tilde b_k \Vert_2^2   ] } = \frac{ \Vert \tilde A_k u \Vert_2^2 }{u^\top \ex [ \tilde A_k^\top \tilde A_k ] u} \leq \frac{ \Vert \tilde A_k \Vert_2^2 }{ \sigma_{\min}( \ex[ \tilde A_k^\top \tilde A_k] )},
    \end{equation}
    where $\sigma_{\min}$ denotes the smallest non-zero singular value of the given matrix. Using \ref{eq:collocA}, we can find the maximum possible value of $\tilde A_k$ on the unit cube is $3$; thus, using the equivalence between the infinity and two norms $\|\tilde A_k\|_2^2 \leq 9N$. So we can set $V = \frac{9N}{\sigma_{\min}( \ex[ \tilde A_k^\top \tilde A_k])}$.
    
    Hence,
    \begin{equation}
    -1 \leq \frac{ \Vert \tilde A_k x - \tilde b_k \Vert_2^2 }{ \ex[ \Vert \tilde A_k x - \tilde b_k \Vert_2^2   ] } - 1\leq V - 1.
    \end{equation}
    Applying \cref{lemma:boundedSE} gives the conclusion.
\end{proof}

It should be emphasized that this result holds true for any right-hand side provided that it is bounded over the domain of the collocation problem. To generalize this result to different right-hand sides one can replace the $9$ in $V$ with the maximum of the right-hand side over the desired domain. Additionally, it should be noted that in practice $V/2$ is far too loose of a bound on to $\sigma$ and should be replaced with estimates from more practical techniques, such as those presented in \cref{subsec:boundary}.

\setcounter{section}{1}
\end{document}